\documentclass[twocolumn,8pt]{article}
\usepackage{times,latexsym,amsmath,amssymb}
\usepackage{graphicx}
\newtheorem{definition}{Definition}[section]
\newtheorem{remark}{Remark}[section]
\newtheorem{lemma}{Lemma}[section]
\newtheorem{theorem}{Theorem}[section]
\newtheorem{proposition}{Proposition}[section]

\newtheorem{example}{Example}[section]

\begin{document}
\global\def\refname{{\normalsize \it References:}}
\baselineskip 12.5pt
%
%
% TITLE, AUTHOR, ABSTRACT, KEYWORDS
%
\title{\LARGE \bf Parameterized Norm and Parameterized Fixed-Point Theorem\\  by Using Fuzzy Soft Set Theory}

\date{}

\author{\hspace*{-10pt}
\begin{minipage}[t]{2in} \normalsize \baselineskip 12.5pt
\centerline{Azadeh Zahedi Khameneh} \centerline{University Putra Malaysia} \centerline{Institute for Mathematical Research} \centerline{43400 UPM Serdang,
Selangor} \centerline{MALAYSIA} \centerline{azadeh503@gmail.com}
\end{minipage} \kern 0in
\begin{minipage}[t]{2in} \normalsize \baselineskip 12.5pt
\centerline{Adem Kili$\c{c}$man} \centerline{University Putra Malaysia} \centerline{Department of  Mathematics} \centerline{43400 UPM Serdang, Selangor}
\centerline{MALAYSIA} \centerline{akilicman@yahoo.com }
\end{minipage}
\begin{minipage}[t]{2in} \normalsize \baselineskip 12.5pt
\centerline{Abdul Razak Salleh} \centerline{University Kebangsaan Malaysia} \centerline{School of  Mathematical Science} \centerline{43600 UKM Bangi,
Selangor} \centerline{MALAYSIA} \centerline{aras@ukm.my}
\end{minipage}
% If you are three authors then you can use three mini--pages
% instead of two. Their horizontal size must be less than 2.7in
% indicated above. It can be e.g. 2.3in. However, you must pay
% attention that you do not exceed the total width of the text.
%
\\ \\ \hspace*{-10pt}
\begin{minipage}[b]{6.9in} \normalsize
\baselineskip 12.5pt {\it Abstract:}
% The text of the abstract follows.
From last decade, when Molodtsov introduced the theory of soft set as a new approach to deal with uncertainties, until now this theory was considered
sharply by a fair number of researchers. Combination of fuzzy set theory and soft set theory, called fuzzy soft set theory, by Maji et.al opened a new way
for researchers whose frame work of study is soft sets and fuzzy sets. Although published papers in this area have considered both application and
theoretical aspects of fuzzy soft set theory, the concept of norm for a fuzzy soft set has not been studied yet. We begin this paper by introducing fuzzy
soft real numbers, which are needed for study fuzzy soft norm, and then continued by considering fuzzy soft norm. Fixed-point theorem is also investigated
for fuzzy soft normed spaces.
\\ [4mm] {\it Key--Words:}
% The key-words follow.
soft set, fuzzy soft set, fuzzy soft topology, fuzzy soft point, fuzzy soft real number, fuzzy soft norm, fuzzy soft continuous map, fuzzy soft fixed-point
theorem
\end{minipage}
\vspace{-10pt}}

\maketitle

\thispagestyle{empty} \pagestyle{empty}
% numbers of pages are supplemented by the editor
%
% THE BEGINNING OF THE TEXT
%
\section{Introduction}
\label{S1} \vspace{-4pt}
 We live in the world of uncertainties, where most of the problems which we face are vague rather than precise. These vague concepts
are due to limited knowledge or incomplete information. During past decades different mathematical theories like probability theory, fuzzy set theory
\cite{Z}, and rough set theory \cite{paw} were introduced to deal with various types of uncertainties. Although a wide range of problems can be solved by
these methods, some difficulties has been remained.\\
Lack of parameterization tools in all previous theories, leaded to introduce soft set theory by Molodtsove \cite{M} in 1999. A soft set is in fact a
set-valued map which gives an approximation description of objects under consideration based on some parameters. Hence the set of all soft sets over a
universal set can be considered as a function space. Since then, Maji et.al \cite{M-B-R,M-R} discussed theoretical aspects and practical applications of
soft sets in decision making problems. For more details about soft set theory see \cite{M-B-R,M-R,A-C,A-F-M,S}. But since almost all the time in real life
situations we have inexact information about considered objects, in 2002, Maji et.al  \cite{m-b-r} studied the combination of fuzzy set theory and soft set
theory and gave a new concept called fuzzy soft set. This new notion expanded the concept of soft set from crisp cases to fuzzy cases and then in
\cite{m-r}, they applied this new method to solve decision making problems. Kharal and Ahmad in \cite{BA-AKH} studied some properties of fuzzy soft set and
then in \cite{kh-a} defined the concept of mapping on fuzzy soft classes. Topological studies of fuzzy soft sets was started by Tanay and Kandemir in
\cite{T-K}. Mahanta and Das \cite{M-D} continued working on fuzzy soft topology and studied separation axioms and connectedness of fuzzy soft topological
spaces. Simisekler and Yuksel \cite{S-Y} modified the definition of fuzzy soft topology and also gave the concept of soft quasi neighborhood of a fuzzy
soft point. Roy and Samanta \cite{R-S} gave the new definition of fuzzy soft topology and proposed the concept of base and subbase for
this new space.\\
Although different aspects of fuzzy soft set theory has been studied by several authors, fuzzy soft number and fuzzy soft norm have not been investigated
yet. In real life situations, we usually use some phrases or sentences like '' distance between A and B is about 10 kilometers'' or '' distance between A
and B is around 20 minutes''. Such vague expressions are used when we deal with quantities which we do not know their values precisely. Moreover, these
inexact values are usually depend on some parameters. For instance, distance between A and B can be changed based on selected measure.
 In addition, it can be presented longer or shorter with regards to chosen path or time of movement. Whereas in classic mathematics distance between two objects in a normed space is
supposed constant, in real life it is not a fixed concept. To
solve such problems a parameterization version of numbers and norm shall be needed. \\
 In this paper, we first recall the definition of fuzzy soft topology which has been established in \cite{R-S}
and then answer to the natural question '' is there any
 relation between point-set topology, fuzzy topology and fuzzy soft topology over a common universe?''.
  We continue our work by introducing the concept of fuzzy
 soft real numbers and then initiate fuzzy soft norm over a set. We also consider the relationship between fuzzy soft norm and fuzzy norm over a common
 set. Finally, we will consider fuzzy soft fixed-point theorem in fuzzy soft normed spaces.
%%%%%%%%%%%%%%%%%%%%%%%%%%%%%%%%%%%%%%%%%%%%%%%%%%%%%%%%%%%%%%%%%%%%%%%%%

\section{Preliminaries}\label{construction} \vspace{-4pt}
Let $X$ be the set of all objects and $E$ be the set of all parameters. Let $A\subseteq E$, and $2^X$ and $I^X$, where $I=[0,1]$, denote the set of all
subsets and fuzzy subsets of $X$, respectively.
 Molodtsov \cite{M} introduced the concept of soft set as follow:
\begin{definition}{\bf:}
(\cite{M}) A pair $(F,A)$ is called a soft set over $X$ if $F$ is a mapping given by $F : A\rightarrow 2^X$ such that for all $e\in A$, $F(e)\subseteq X$.
$\label{def2.1}$
\end{definition}
The soft set $(F,A)$ can be denoted by $F_A$.\\

 Maji et.al \cite{m-b-r} gave the concept of fuzzy soft set as below:
\begin{definition}{\bf:}
(\cite{m-b-r}) A pair $(f,A)$ is called a fuzzy soft set over $X$, F.S set briefly, if $f$ is a mapping given by $f : A\rightarrow  I^X$. So $\forall e\in
A$, $f(e)$ is a fuzzy subset of $X$, with membership function
$$f_e: X\rightarrow [0,1]$$  $\label{def2.2}$
\end{definition}
In fact, the membership function $f_e$ indicates degree of belongingness of each element of $X$ in $f(e)$ or shows how much each member of $X$ has the
parameter $e\in E$.
 So the soft set $(f,A)$ can be represented by the set of triplet
ordered $$\{(e,x,f_e(x)): e\in A , x\in X , f(e)\in I^X\}$$ where $f_e(x)$ is degree of
membership $x$ in fuzzy set $f(e)$.\\
 The fuzzy soft set $(f,A)$ can be denoted by $f_A$.

\begin{example}
In geography science, definition of a forest is expressed based on some parameters such as domain and density of vegetation, type of vegetation, soil
types, amount of annual rain, and species of plant and animal. Let $X=\{A,B,C\}$ be a set of regions under consideration. Suppose that
$E$=$\{e_1,e_2,e_3,e_4\}$ is a set of parameters where $e_i$ $(i=1,\ldots, 4)$ stand for the parameters: wide, dense vegetation, rainy, variety of plant
and animal, respectively.
 The available information of these regions can be presented by the fuzzy soft set $f_E$ as below:\\
  $$f_E= \{ (e_1,\{\frac{0.8}{A} , \frac{0.3}{B} , \frac{0.5}{C}\}) , (e_2,\{\frac{0.1}{A} , \frac{0.5}{B} , \frac{0.7}{C}\}) ,$$
  $$(e_3,\{\frac{0.2}{A} , \frac{0.3}{B} , \frac{0.8}{C}\})
   , (e_4,\{\frac{0.1}{A} , \frac{0.3}{B} , \frac{0.5}{C}\}) \}$$
 which means that $A$ is the largest zone among these three areas, while it has the lowest amount of vegetation, annual rain, and
diversity of animal and plant species.
     $\label{ex2.1}$
\end{example}

\begin{definition}{\bf:}
(\cite{m-b-r},\cite{kh-a}) (Rules of fuzzy soft set) For two fuzzy soft sets $f_A$ and $g_B$ over the common universe $X$ with respect to parameter set $E$
where $A,B \subseteq E$ we have,
\begin{description}
    \item[i.] $f_A$ is a fuzzy soft subset of $g_B$ shown by
$f_A\tilde\leq g_B$ if:
\begin{enumerate}
\item $A\subseteq B$, \item For all $e$ in $A$, $f_e (x)\leq g_e (x)$, $\forall x\in X$.
\end{enumerate}
    \item[ii.] $f_A = g_B$ if $f_A\tilde\leq g_B$ and
    $g_B\tilde\leq f_A$.
    \item[iii.]  The complement of fuzzy soft set $f_A$ is denoted by $f_A^{c}$
 where $f^{c}: A\rightarrow I^X$
 and $f^{c}(e)$ is the complement of fuzzy set $f(e)$, i.e $f^{c}(e)=(f(e))^c$, with
 membership function $f^{c}_e=1-f_e$, $\forall e\in A$.
    \item[iv.] $f_A=\Phi_A$ (null fuzzy soft set with respect to $A$), if for each $e\in A$,
$f_e(x)=0$, $\forall x\in X$.
    \item[v.] $f_A=\tilde X_A$ (absolute fuzzy soft
set with respect to $A$) if if for each $e\in A$, $f_e(x)=1$ , $\forall x\in X$.\\
 If $A=E$, the null and absolute fuzzy soft set is denoted by $\Phi$ and $\tilde
 X$, respectively.
    \item[vi.] The union of two fuzzy soft sets $f_A$ and $g_B$, denoted by $f_A\tilde\vee g_B$, is the fuzzy
 soft set $(f\vee g)_C$ where $C=A\cup B$ and $\forall e\in C$, we have $(f\vee g)(e)=f(e)\vee g(e)$ where\\
 $(f\vee g)_e(x)=$
$ \left\{\begin{array}{ll}
 f_e(x) & \mbox{if $e\in A-B$}\\
 g_e(x) & \mbox{if $e\in B-A$}\\
 \max\{f_e(x) , g_e(x)\} & \mbox{if $e\in A\cap B$}
 \end{array}
 \right.$\\
 \ for all $x\in X$.
    \item[vii.] The intersection of two fuzzy soft sets $f_A$ and $g_B$, denoted by $f_A\tilde\wedge g_B$, is the fuzzy soft set $(f\wedge g)_C$
     where $C=A\cap B$ and $\forall e\in C$, we have $(f\wedge
g)(e)=f(e)\wedge g(e)$ where $(f\wedge g)_e(x)=\min\{f_e(x),g_e(x)\}$ for all $x\in X$.
\end{description}
 $\label{def2.3}$
\end{definition}
In (vii), $A\cap B$ must be nonempty to avoid the degenerate case.\\
 Note that during this paper, $X$ and $E$ are used to show the universal sets of objects and
parameters, respectively. $f_A$ denotes a fuzzy soft set over $X$ where $A\subseteq E$ and $X_E$ denotes the set of all fuzzy soft sets over $X$ with
regards to parameter set $E$.

\begin{proposition}
(\cite{A-K-S}) (De Morgan Laws) Let $f_E$ and $g_E$ be F.S sets over $X$ with respect to $E$. Then we have
\begin{enumerate}
    \item $[f_E\tilde\vee g_E]^c=f_E^c\tilde\wedge g_E^c$
    \item $[f_E\tilde\wedge g_E]^c=f_E^c\tilde\vee g_E^c$
\end{enumerate}
$\label{pr3.1}$
\end{proposition}
\textbf{Proof.} See \cite{A-K-S}. $\Box$

\begin{definition}{\bf:}
(\cite{kh-a}) Let $X$ and $Y$ be universal sets, and $E$ and $E'$ be corresponding parameter sets, respectively.
    Suppose that $f_A$ be a F.S set on $X$ and $g_B$ be a F.S set on $Y$. Let
     $u:X\rightarrow Y$ and $p: E\rightarrow E'$ be ordinary functions.
\begin{description}
    \item[i.] The map $h_{up}: X_E \rightarrow
Y_{E'}$ is called a F.S map from $X$ to $Y$ mapping $f_A$ to F.S set $h_{up} (f_A)$ and for all $y\in Y$ and $e'\in p(E)$ it is defined as below:
 $$[h_{up} (f)]_{e'}(y) =
 \sup_{x\in u^{-1}(y)}[\sup_{e\in p^{-1}(e')\cap A} f(e)] (x)$$ if $p^{-1}(e')\cap A\neq\emptyset,u^{-1}(y)\neq\emptyset$\\ and otherwise $[h_{up}
 (f_A)]_{e'}(y)=0$.
    \item[ii.] Let $h_{up}: X_E \rightarrow
Y_{E'}$ be a F.S map from $X$ to $Y$. Then the inverse image of F.S set $g_B$, denoting by $h^{-1}_{up} (g_B)$,  is a F.S set on $X$ and for all $x\in X$
and $e\in E$ is defined as below:
 $$[h^{-1}_{up} (g)]_{e}(x)= \left\{\begin{array}{ll}
 g_{p(e)}(u(x)) & \mbox{if $p(e)\in B$}\\
 0 & \mbox{otherwise}\\
 \end{array}
 \right.$$
 \end{description}
 $\label{def7.0}$
\end{definition}

\begin{theorem}{\bf:}
(\cite{kh-a})  Let $h_{up}: X_E \rightarrow Y_{E'}$ be a F.S map where $u:X\rightarrow Y$ and $p: E\rightarrow E'$ are ordinary functions. For fuzzy soft
set $f_A$ and family of fuzzy soft sets $(f_A)_i$ in $X$ and for fuzzy soft set $g_B$ and
    family of fuzzy soft sets $(g_B)_i$ on $Y$ we have
\begin{enumerate}
    \item $h_{up}(\Phi)=\Phi$ and $h_{up}^{-1}(\Phi)=\Phi$
    \item $h_{up}(\tilde X)\tilde\leq\tilde Y$ and $h^{-1}_{up}(\tilde Y)=\tilde X$
    \item $h_{up}[\tilde\bigvee_i (f_A)_i]=\tilde\bigvee_i h_{up}(f_A)_i$ and\\
      $h^{-1}_{up}[\tilde\bigvee_i (g_B)_i]=\tilde\bigvee_i h^{-1}_{up}(g_B)_i$
    \item $h_{up}[\tilde\bigwedge_i (f_A)_i]\tilde\leq\tilde\bigwedge_i h_{up}(f_A)_i$ and\\
      $h^{-1}_{up}[\tilde\bigwedge_i (g_B)_i]=\tilde\bigwedge_i h^{-1}_{up}(g_B)_i$
    \item $[h_{up}(f_A)]^c\tilde\leq h_{up}(f_A)^c$ and $[h^{-1}_{up}(g_B)]^c = h^{-1}_{up}(g_B)^c$
    \end{enumerate}
 $\label{th7.0}$
\end{theorem}

\textbf{Proof.} See \cite{kh-a}. $\Box$\\

Here we apply Definition $\ref{def7.0}$ and Fuzzy Extension Principle for an ordinary function (see \cite{ch}) to give a parameterized extension of an
ordinary function to a F.S map.
\begin{definition}{\bf:}
Let $u: X\rightarrow Y$ be a function. Let $f_A$ be a F.S set in $X$. Then the image of of $f_A$ is a F.S set in $Y$ whose membership function is given by:
 $$[u (f)]_{e'}(y) =  \sup_{x\in u^{-1}(y)}[\sup_{e\in p^{-1}(e')\cap A} f(e)] (x)$$
 if $p^{-1}(e')\cap A\neq\emptyset,u^{-1}(y)\neq\emptyset$ and otherwise $0$. This is shown in the following Diagram\\
  \begin{eqnarray*}
E & \longrightarrow^p & E'\\
\downarrow^f & & \downarrow^{u(f)}\\
X & \longrightarrow^u & Y\\
f(e)^{\searrow} & & \swarrow_{f(e)ou^{-1}}\\
 &[0,1]& \\
  \end{eqnarray*}
 Hence $(uf)(e')=f(p^{-1}(e'))o u^{-1}$.\\
    Let $g_{B}$ be a F.S set in $Y$. The inverse image of $g_{B}$ is a F.S set in $X$ whose membership function is given by:
 $$[u^{-1} (g)]_{e}(x)= \left\{\begin{array}{ll}
 g_{p(e)}(u(x)) & \mbox{if $p(e)\in B$}\\
 0 & \mbox{otherwise}\\
 \end{array}
 \right.$$
   This is shown in the following Diagram\\
   \begin{eqnarray*}
E & \longrightarrow^p & E'\\
 {u^{-1}g}{\downarrow} & & \downarrow^{g}\\
X & \longrightarrow^u & Y\\
{g(e')o u} {\searrow} & & \swarrow_{g(e')}\\
 &[0,1]& \\
  \end{eqnarray*}
  So $g(p(e))o u=(u^{-1}g)(e)$.
  $\label{def7.0.1}$
\end{definition}
 Note that, If $p: E \rightarrow E'$ is a bijective mapping, then the former definition can be written as below:
 $$[u (f)]_{e}(y) = \left\{\begin{array}{ll}
  \sup_{x\in u^{-1}(y)} f_e(x) & \mbox{$e=p^{-1}(e')\in A,u^{-1}(y)\neq\emptyset$}\\
 0 & \mbox{otherwise}
 \end{array}
 \right.$$
 since there exists only one element in the parameter set $E$ such that $p(e)=e'$. If $A=E$, then we have
 $$[u (f)]_{e}(y) = \left\{\begin{array}{ll}
  \sup_{x\in u^{-1}(y)} f_e(x) & \mbox{$u^{-1}(y)\neq\emptyset$}\\
 0 & \mbox{otherwise}
 \end{array}
 \right.$$
 %%%%%%%%%%%%%%%%%%%%%%%%%%%%%%%%%%%%%%%%%%%%%%%%%%%%%%%%%%%%%%%%%%%%%%%%%
\section{Fuzzy Soft Topology and its Relation with Point-Set and Fuzzy Topologies} \vspace{-4pt}
 Topological studies of fuzzy soft sets was began by Tanay and Kandemir in \cite{T-K}. The concept of fuzzy soft topology
introduced by them is in fact, a topological structure over a fuzzy soft set. So it can be seen as a collection of fuzzy soft subsets of an arbitrary fuzzy
soft set. This means that parameter set is not fixed everywhere (see Definition $\ref{def2.3}$ part (i)). But since in this case De Morgan law's are not
hold in general (see \cite{A-K-S}),
   in \cite{R-S}, Roy and Samanta initiated the concept of fuzzy soft topology over a universal set
where the parameter set is supposed fixed all over the universe. Here we recall definition of fuzzy soft topology introduced in \cite{R-S}, and consider
the relationship between this new topology with two previous topologies, point-set topology and fuzzy topology, over a common universe.

\begin{definition}{\bf:}
(\cite{R-S}) A fuzzy soft topology over $X$ denoted by $\tau$, is a collection of fuzzy soft subsets of
 $X$
such that:
\begin{description}
    \item[i.] $\tilde X$ and
$\Phi\in\tau$.
    \item[ii.] The union of any number of fuzzy soft sets
in $\tau$ belongs to $\tau$.
    \item[iii.] The intersection of any two fuzzy soft sets in $\tau$ belongs
to $\tau$.
\end{description}
$\label{def4.1.1}$
\end{definition}
 The triplet $(X,E,\tau)$ is called a fuzzy soft
topological space, F.S topological space in brief, and
 each element of $\tau$ is called a fuzzy soft open set, say F.S open set, in $X$.
 The complement of a F.S open set is called F.S closed set.

\begin{example}{\bf:}
Let $X$ be a universal set including objects under consideration and $E$ be the set of parameters. Then the family of all F.S sets over $X$, denoting by
$X_E$, forms a F.S topology over $X$ which is called discrete F.S topology, while indiscrete or trivial F.S topology on $X$ contains only $\Phi ,\tilde X$.
$\label{ex4.1.1}$
\end{example}

\begin{example}{\bf:}
Let $(X,\tau)$ be a topological space. Then the family
 $$\{V_E : V\in \tau , E=(0,1]\}$$
  forms a F.S topology over $X$ denoting by
$\tau_{F.S}$ where $V_E$ is a F.S set over $X$ with respect to $\tau$-open set $V$ defined as below:
 $$ V: E=(0,1]\rightarrow I^{X}$$
  where for each $\alpha\in (0,1]$, $V(\alpha)$ is defined by characteristic function of
  $\tau$-open set $V$ i.e., $V(\alpha)=\chi_{V}$. So\\
$$V_\alpha(x)= \left\{\begin{array}{ll}
 1 & \mbox{ $x\in V$}\\
 0 & \mbox{$x\notin V$}\\
 \end{array}
 \right.$$
$\label{ex4.2.1}$
\end{example}

\begin{example}{\bf:}
Let $(X,\tau)$ be a topological space and $E$ be a parameter set. The collection $\tau_{F.S}=\{f; f: E \rightarrow I^X: \forall e\in E,
(f_e)^{-1}(0,1]\in\tau \}$ is a fuzzy soft topology over $X$ where for each $e\in E$, $(f_e)^{-1}(0,1]$ denotes the support of fuzzy set $f(e)$. Moreover
$\tau\subset \tau_{FS}$. $\label{ex5.2}$
\end{example}

\begin{example}{\bf:}
Let $(X,\gamma)$ be a fuzzy topological space. Let $\mu$ be a fuzzy subset of $X$ and $\chi$ denotes the characteristic function.
 We define the characteristic
function of $\alpha$-cut sets of $\mu$, denoting by $\chi_{\mu_{\alpha}}$, as below:
$$ \chi_{\mu_{\alpha}}: X\rightarrow [0,1]$$
$$(\chi_{\mu_{\alpha}})(x)= \left\{\begin{array}{ll}
 1 & \mbox{$\mu(x)\geq \alpha$}\\
 0 & \mbox{$\mu(x) < \alpha$}\\
 \end{array}
 \right.$$
Now we construct the the F.S map ${\chi_{\mu}}_E$ by using $\chi_{\mu_{\alpha}}$'s as below:
$$ \chi_{\mu}: E=(0,1]\rightarrow I^{X}$$
$$(\chi_{\mu})(\alpha)=\chi_{\mu_{\alpha}}$$
 for each $\alpha\in E$.
Then the collection $$\tau_{F.S}=\{\chi_{\mu_{E}}: \chi_{\mu}: E \rightarrow I^X; \forall \alpha\in E, (\chi_{\mu})(\alpha)=\chi_{\mu_{\alpha}} \}$$ is a
fuzzy soft topology over $X$. $\label{ex5.3}$
\end{example}

\begin{theorem}{\bf:}
(\cite{A-K}) Let $(X,E,\tau)$ be a fuzzy soft topological space. Then corresponding to each $e\in E$, $\tau_e=\{f(e) : f_E\in\tau , e\in E\}$ forms a fuzzy
topology over $X$ (in the sense of Chang \cite{ch}). $\label{th5.1}$
\end{theorem}
\textbf{Proof.}
\begin{description}
    \item[i.] $\Phi , \tilde X\in\tau \Rightarrow \underline{0},\underline{1}\in\tau_e$
    where $\underline{0}$ denotes the empty fuzzy set and $\underline{1}$ shows the fuzzification of $X$, i.e.\\
  $\underline{0}=\chi_{\emptyset}$,
  $\underline{1}=\chi_{X}$,
  and $\chi$ denotes the characteristic function.
    \item[ii.] Let $\Lambda$ be an index set and $\{f^\lambda_E\}_{\lambda\in\Lambda}\subseteq\tau$. Let $e\in E$, then for all $\lambda\in\Lambda$,
    $f^\lambda(e)\in\tau_e$. Since $\tau$ is a F.S topology over $X$, then by applying Definition $\ref{def2.3}$ part (vi) we have
$\tilde\bigvee_{\lambda\in\Lambda}f^\lambda_E\in\tau \Rightarrow (\bigvee_{\lambda\in\Lambda}f^\lambda)_E\in\tau  \Rightarrow
(\bigvee_{\lambda\in\Lambda}f^\lambda) (e)\in\tau_e$.
    \item[iii.] Suppose $f^1_E , f^2_E\in \tau$ and $e\in E$. Then $f^{1}(e) ,
f^{2}(e)\in\tau_e$. Since $\tau$ is a F.S topology, then then by applying Definition $\ref{def2.3}$ part (vii) we have
 $f^1_E\tilde\wedge f^2_E\in\tau \Rightarrow (f^1\wedge f^2)_E\in\tau$ $\Rightarrow
f^{1}(e)\wedge f^{2}(e)\in\tau_e$.\\
So $\tau_e$ is a fuzzy topology over $X$.
\end{description}
$\Box$\\

Now we define the concept of fuzzy soft point over set $X$. In the literature, the concept of fuzzy soft point was introduced as bellows:
\begin{enumerate}
    \item (\cite{M-D}) A fuzzy soft set $f_E$ over $X$ is called a fuzzy soft point if for the element
    $e^*\in E$ we have
     $$f_e(x) = \left\{\begin{array}{ll}
\lambda_x & \mbox{if $e=e^*$}\\
0 & \mbox{otherwise}
  \end{array}
 \right.$$
 $\forall e\in E$, $\forall x\in X$ where $\lambda\in (0,1]$.
    \item (\cite{S-Y})  A fuzzy soft set $f_E$ over $X$ is called a fuzzy soft point if for the element
    $x^*\in X$ we have
     $$f_e(x) = \left\{\begin{array}{ll}
\lambda & \mbox{if $x=x^*$}\\
0 & \mbox{otherwise}
  \end{array}
 \right.$$
 $\forall e\in E$, $\forall x\in X$ where $\lambda\in (0,1]$.
\end{enumerate}
But these definitions are not free of difficulties. First one investigates the cases which are related to only one parameter. Second one, although is more
general than the former, considers specific cases in which the membership degree is supposed fixed for all parameters. In fact, it can be seen as a
parameterized version of a fuzzy point, whereas fuzzy soft set theory is a new method to give a fuzzy extension of soft set theory.\\
 In classic set theory, each member of a non-empty set is defined as an object which has the same property of the set.\\
  On the other hand, fuzzy soft set theory is a method used to represent the imprecise information about a set based on some parameters.
 So to define the concept of fuzzy soft point, we restrict the universal set $X$ to the single set $\{x\}\subset X$ and define the fuzzy soft point $x$ as a fuzzy description of
 element $x\in X$ with regards to some parameters.

\begin{definition}{\bf:}
(\cite{A-K,A-K-S}) \begin{enumerate} \item Let $x^{\lambda_e}$ be a fuzzy point with support $x\in X$ and membership degree $\lambda_e\in(0,1]$ (see
\cite{pu-liu} ). The fuzzy soft set ${Px}_E$ is called fuzzy soft point, say F.S point, whenever $Px: E\rightarrow I^X$ is a map such that for each $e\in
E$ and $\forall z\in X$
$$(Px)_e(z) =
\left\{\begin{array}{ll}
\lambda_e & \mbox{if $z=x$}\\
0 & \mbox{otherwise}
  \end{array}
 \right.$$
So for each $e\in E$, $(Px)(e)=x^{\lambda_e}$, or $(Px)(e)=\lambda_e \chi_{\{x\}}$ where $\chi_{\{x\}}$ is the characteristic function of $\{x\}$.
    In other words, the F.S point ${Px}_E$, is a
    fuzzy description of $x\in X$ based on parameter
  set $E$.\\
    If $\lambda_e=1$ for all $e\in E$, we call ${Px}_E$, crisp F.S
    point.
        \item The F.S point ${Px}_E$  belongs to F.S
    set $f_E$ denoting by ${Px}_E\tilde\in f_E$,
    whenever for all $e\in E$ we have $0<\lambda_e\leq f_e(x)$.
    \item The restriction of F.S point ${Px}_E$ to an element $e\in E$, denoting by ${Px}_E|_e$, is
called fuzzy soft single point over $E$, say F.S single point, whenever for all $\alpha\in E$ and $\forall z\in X$
$$(Px|_e)_\alpha(z) = \left\{\begin{array}{ll}
  (Px)_e(x)=\lambda_e & \mbox{if $\alpha=e$ , $z=x$}\\
0 & \mbox{otherwise}
  \end{array}
 \right.$$
The F.S single point ${Px}_E|_e$  belongs to F.S
    set $f_E$ denoting by ${Px}_E|_e\tilde\in f_E$,
    whenever for $e\in E$ we have $0<\lambda_e\leq f_e(x)$.
  \end{enumerate} $\label{def4.2.1}$
\end{definition}
It is clear that Definition $\ref{def4.2.1}$ is an extension of both crisp and fuzzy points. In addition, the concept of fuzzy soft point
 introduced in the literature \cite{M-D,S-Y} can be seen as a specific case of our definition.\\

  \textbf{New Notation.}
\begin{itemize}
    \item The F.S point
${Px}_E$  will be denoted by
     $\tilde x_E$ where for all $e\in E$, $(Px)(e)=\tilde x(e)=x^{\lambda_e}$ i.e.,
     the image of each parameter under map $Px$ is a fuzzy point. Consequently  $(\tilde x_E)^c$
    can be applied to show the complement of F.S point $\tilde x_E$ such that for all $e\in E$ and $\forall z\in X$, we have
    $$(\tilde x_e^c)(z)=\left\{\begin{array}{ll}
    1-\lambda_e & \mbox{$z=x$}\\
    1 & \mbox{$z\neq x$}
     \end{array}
 \right.$$
    \item The crisp F.S point ${Px}_E$ will be denoted by $x^1_E$ or $\bar{x}_E$.
    \item The F.S single point ${Px}_E|_e$  denoting by $\tilde x_e$.
\end{itemize}
From now, notation $\tilde x_E$ means F.S point $x$ in $X$ where $\tilde x(e)=x^{\lambda_e},\forall e\in E$. \\
 Note that crisp point $x$ and fuzzy point $x^\lambda$  can be viewed as a fuzzy soft points $\bar{x}_{(0,1]}$ and $\tilde x_{(0,1]}$ where
$\forall \alpha\in E=(0,1]$, $\bar x(\alpha)=x^1$ and $\tilde x(\alpha)=x^\lambda$, respectively %(see Fig. $\ref{Fig.1}$).
%\begin{figure}
% \includegraphics[width=3.25 in]{Doc4}
% \caption{Crisp F.S point $\bar{x}_{(0,1]}$}\label{Fig.1}
%\end{figure}

\begin{definition}
(\cite{A-K,A-K-S}) Let $(X,E,\tau)$ be a fuzzy soft topological space. Two fuzzy soft points $\tilde x_E$ and $\tilde y_E$ where for all $e\in E$, $\tilde
x(e)=x^{\lambda_e}$ and $\tilde y(e)=y^{\gamma_e}$  are said to be
\begin{itemize}
    \item  different if and only if
\begin{enumerate}
    \item  $x\neq y$ or
    \item Whenever $x=y$, we have  $\lambda_e \neq
    \gamma_e$ for some $e\in E$.
\end{enumerate}
    \item distinct if and only if
 $\tilde x_E \tilde \wedge \tilde y_E=\Phi$.
\end{itemize}
$\label{def4.1}$
\end{definition}

\begin{definition}{\bf:}
(\cite{A-K}) Let $\tilde x_E$ be a F.S point over $X$. Fuzzy soft set $g_E$ is called  fuzzy soft neighborhood, F.S - N in brief, of F.S point $\tilde x_E$
whenever there exists a F.S open set $f_E$ such that $\tilde x_E\tilde\in f_E\tilde\leq g_E$. $\label{def4.2.2}$
\end{definition}

\begin{example}{\bf:}
Take the set of all real numbers $\mathbb{R}$ with usual topology $\tau_u$. Let $V=(a,b)\subseteq \mathbb{R}$ be an open neighborhood of $x\in \mathbb{R}$.
 Define fuzzy soft set $V_E$ as the following:
 \begin{eqnarray*}
  V: (0,1] &\rightarrow& I^{\mathbb{R}}\\
  \alpha &\mapsto& V(\alpha)=\chi_{V}
  \end{eqnarray*}
Then $V_E$ is an F.S-open-N of F.S point $\tilde x_E$ in F.S topological space $(X,(0,1],\tau_{F.S})$ where $\tau_{F.S}$ is the F.S topology mentioned in
Example $\ref{ex4.2.1}$. In fact, $V_E$ is a parametrization version of  open interval $V$ in $\mathbb{R}$ %(see Fig. $\ref{Fig.2}$).
\end{example}

%\begin{figure}
% \includegraphics[width=3.25 in]{Doc2}
% \caption{Parametrization of open interval (a,b)}\label{Fig.2}
%\end{figure}

\begin{definition}{\bf:}
(\cite{A-K})
\begin{enumerate}
    \item Let $\tilde x_E$ be a F.S point in $X$.
    We say that $\tilde x_E$ is soft quasi-coincident with F.S set $f_E$, denoting by $\tilde x_E \tilde{q} f_E$, if there exists
$e\in E$, such that $\lambda_e + f_e(x)>1$.
    \item The fuzzy soft set $g_E$ is called soft quasi-coincident with fuzzy soft set $f_E$ at $x\in X$, denoting by $g_E\tilde{q} f_E$,
    if there exist $e\in E$ such that $g_e(x) + f_e(x)>1$. If not we say that $g_E$ is not soft quasi-coincident with $f_E$ denoting by
   $f_E\neg{\tilde {q}}g_E$.
    \item The fuzzy soft single point $\tilde x_e$ is called soft quasi-coincident with fuzzy soft set $f_E$ at $x$, denoting by $\tilde x_e \tilde{q}
f_E$, if for $e\in E$ we have $\lambda_e + f_e(x)>1$.
\end{enumerate}
   $\label{def4.2.3}$
\end{definition}

\begin{definition}{\bf:}
(\cite{A-K}) The fuzzy soft set $g_E$ is called a soft quasi neighborhood of F.S point $\tilde x_E$, say soft Q  - N, if there exists the F.S open subset
$f_E$, such that $\tilde x_E\tilde{q} f_E \tilde\leq g_E$. $\label{def4.2.4}$
\end{definition}

Let $\tilde x_E$ be a F.S point in $X$. Let $\Lambda$ be an index set and $(f_E)_{\alpha\in \Lambda}$ be a family of fuzzy soft sets over $X$ with respect
to parameter set $E$.

\begin{proposition}
If $\tilde x_E\tilde\in\tilde\bigwedge_{\alpha\in\Lambda} (f_E)_\alpha \Rightarrow \forall\alpha\in\Lambda: \tilde x_E \tilde\in (f_E)_\alpha$.
\end{proposition}
\textbf{Proof.} See \cite{A-K}.

\begin{proposition}
If $\tilde x_E \tilde q \tilde\bigvee_{\alpha\in\Lambda} (f_E)_\alpha \Rightarrow \exists \alpha\in\Lambda: \tilde x_E \tilde q (f_E)_\alpha$.
\end{proposition}
\textbf{Proof.} See \cite{A-K}.

\begin{definition}
(\cite{A-K,A-K-S}) Let $(X,E,\tau)$ be a fuzzy soft topological space. We say that $X$ is a fuzzy soft
\begin{enumerate}
    \item $T_0$, say F.S $T_0$, if and only if for every two distinct F.S points in $X$,
 at least one of them has a F.S
open - N which is  not intersection with the other.
    \item $T_1$ space, say F.S $T_1$, if and only if for every two distinct F.S
points in $X$ such as
 $\tilde x_E$ and $\tilde y_E$,  there exist two F.S
open - N of $\tilde x_E , \tilde y_E$ like $f_E$ and $g_E$ respectively, such that $\tilde y_E \tilde\wedge f_E=\Phi$, and $\tilde x_E\tilde\wedge
g_E=\Phi$.
    \item Hausdorff space, say F.S $T_2$ or F.S Hausdorff, if and only if for
every two distinct F.S points in $X$ such as $\tilde x_E$ and $ \tilde y_E$,  there exist two F.S open - N like $f_E$ and $g_E$, such that $f_E\tilde
\wedge g_E=\Phi$.
\end{enumerate}
 $\label{def4.2}$
\end{definition}

Now we introduce a way to construct a fuzzy soft topology over a set  associate with the given point-set topology or fuzzy topology over
it.\\
We know that for each ordinary set like $A$
 if $A$ is an uncountable set, then we have $A\sim \mathbb{R}\sim (0,1]$,
 and if $A$ is a finite or countable set, then  $A\sim \mathbb{N}_k$ or $A\sim \mathbb{N}$, respectively where $k$ is the number of elements of $A$.\\
 Now consider the parameter set $E$. In the following remark we discuss how the parameter set $E$ can be replaced with a suitable subset of $(0,1]$.
\begin{remark}{\bf:}
(\cite{A-K}) Let $X$ be the set of objects under consideration, and $E$ be the set of parameters.
\begin{enumerate}
    \item If $E$ is an uncountable set, then we replace $E$ with
    $(0,1]$. So associate with each $e\in E$ we have an $\alpha\in (0,1]$.
    \item If $E$ is an infinite but countable set, we replace $E$ with
    the set of all rational numbers in $(0,1]$. Thus in this case, each $e_i\in E$ can be replaced with
    $\frac{m}{n}\in (0,1]$ such that $m=1,2,\ldots$ , $n=1,2,\ldots$,
     and $m \leq n$.
    \item If $E$ is a finite set where $|E|=k$, we replace $E$ with
    $\{\frac{1}{k},\ldots,\frac{k-1}{k},1\}$.
     So corresponding to each $e_i\in E$, we have $\frac{i}{k}\in (0,1]$, where $i=1,\ldots,k-1,k$.
      $\label{rem5.1}$
\end{enumerate}
\end{remark}
From now, $(0,1]_E$ is used to show the subset of $(0,1]$ corresponding to $E$.\\
 Now we are ready to introduce a fuzzy soft topology on $X$ associate with
the initial topology on $X$. Put $\mathcal{T}(X)$ the set of all topologies on $X$, $\mathcal{T}_{F}(X)$ the set of all fuzzy topologies on $X$, and
$\mathcal{T}_{F.S}(X)$ the set of all fuzzy soft topologies on $X$ and consider the following mappings
$$i'_e: \mathcal{T}_{F.S}(X) \rightarrow \mathcal{T}_{F}(X) , i_t: \mathcal{T}_{F}(X) \rightarrow \mathcal{T}(X)$$
$$w: \mathcal{T}(X) \rightarrow \mathcal{T}_{F}(X) , w':  \mathcal{T}_{F}(X) \rightarrow \mathcal{T}_{F.S}(X)$$

\begin{definition}{\bf:}
(\cite{A-K}) \begin{enumerate}
    \item Let $(X,E,\delta)$ be a fuzzy soft topological
space, then corresponding to each parameter $e\in E$ we can define a fuzzy topology on $X$ by
 $$i'_e(\delta)=\{f(e) : f_E\in\delta\}$$
 (see Theorem $\ref{th5.1}$) and a topology over $X$ by
  $$i_t(i'_e(\delta))=\{(f(e))^{-1}(t,1] : t\in [0,1), f_E\in\delta\}$$
 as a topology over $X$ where $t\in (0,1]_E$.
    \item Suppose $(X,\tau)$ be a
topological space and $E$ be the set of parameters. We can define a fuzzy soft topology on $X$ with regards to $\tau$ by
 $$w'(w(\tau))=\{f : E \rightarrow I^X ;
f(e)\in w(\tau) ;  \forall e\in E \}$$
 where
  $$w(\tau)=\{\mu : X \rightarrow [0,1] : \mu^{-1}(t,1]\in\tau, \forall t\in[0,1)\}$$
 is a fuzzy topology over $X$ (see \cite{L}).
    \end{enumerate} $\label{def5.1}$
\end{definition}
It is easy to check that $w'(w(\tau))$ and $i_t(i'_e(\delta))$ are indeed  F.S topology  and  point-set topology over $X$, respectively (see \cite{A-K}).
%%%%%%%%%%%%%%%%%%%%%%%%%%%%%%%%%%%%%%%%%%%%%%%%%%%%%%%%%%%%%%%%%%%%%%%%%%%%%%%%%%%%%%%%%%%%%%%%%%%%%
\section{Fuzzy Soft Real Number}
A fuzzy soft number is a fuzzy soft set on real numbers $\mathbb{R}$ associated with some parameters. This notion is a mathematical tool to represent terms
such as '' about $m$ with respect to $n$''. Before giving the concept of fuzzy soft real number, we recall some definitions connected with this subject
which are considered in the literature (see \cite{Z,D-H-L,M-M,K-S}).

\subsection{Fuzzy Real number}
\begin{definition}{\bf:}
(\cite{M-M,K-S}) The fuzzy real number $\mu$ is a fuzzy subset of real number set $\mathbb{R}$ which satisfies the below conditions:
 \begin{itemize}
    \item $\mu$ is normal, i.e., $\exists r\in \mathbb{R}; \mu(r)=1$.
    \item $\mu$ is convex, i.e., $\forall \alpha\in (0,1]$, $\alpha$-level's of $\mu$ are convex sets in $\mathbb{R}$.
    \item $\mu$ is upper semi-continuous, i.e., for all $t\in (0,1]$ and $\varepsilon>0$, $\mu^{-1}([0,t))$ is an open set in $\mathbb{R}$ when topology
    over $\mathbb{R}$ is supposed the usual topology $\tau_u$.
 \end{itemize}
  $\label{7.0}$
\end{definition}
 The set of all fuzzy real numbers is denoted by $\mathcal{F}(\mathbb{R})$. It can be proved that $\alpha$- level sets of an upper semi-continuous convex normal
fuzzy set for each $\alpha\in (0,1]$ is a closed interval in $\mathbb{R}$. So we have
 $$[\mu]_\alpha=\{t\in \mathbb{R}: \mu(t) \geq \alpha\}=[\mu^1_\alpha ,\mu^1_\alpha]$$
  where $\mu^1_\alpha , \mu^1_\alpha\in \mathbb{R}$.\\
A fuzzy number $\mu$ is called non-negative if $\mu(t)=0$ for $t<0$. The set of all non-negative fuzzy real numbers is denoted by
$\mathcal{F}(\mathbb{R^*})$.

\begin{definition}{\bf:}
(\cite{M-M,K-S}) The operations $\oplus, \ominus, \otimes, \oslash$ are defined on $\mathcal{F}(\mathbb{R}) \times \mathcal{F}(\mathbb{R})$ as below:\\
     $(\mu\oplus\delta)(t) = \sup_{s} \min\{\mu(s),\delta(t-s)\}$, $t,s\in \mathbb{R}$\\
     $(\mu \ominus\delta)(t) = \sup_{s} \min\{\mu(s),\delta(s-t)\}$, $t,s\in \mathbb{R}$\\
     $(\mu\otimes\delta)(t) = \sup_{s} \min\{\mu(s),\delta(t/s)\}$, $t,s\in \mathbb{R}$, $s\neq 0$\\
    $(\mu\oslash\delta)(t) = \sup_{s} \min\{\mu(st),\delta(s)\}$, $t,s\in \mathbb{R}$\\
   $\label{7.1}$
\end{definition}

\begin{definition}
  (\cite{M-M,K-S}) We can also define fuzzy arithmetic operations by $\alpha$-level sets as below:\\
$[\mu\oplus\delta]_\alpha = [\mu^1_\alpha+\delta^1_\alpha , \mu^2_\alpha+\delta^2_\alpha]$\\
    $[\mu \ominus\delta]_\alpha = [\min\{\mu^1_\alpha-\delta^1_\alpha,\mu^2_\alpha-\delta^2_\alpha\} ,
     \max\{\mu^1_\alpha-\delta^1_\alpha,\mu^2_\alpha-\delta^2_\alpha\}]$\\
   $[\mu\otimes\delta]_\alpha = [\min\{\mu^1_\alpha.\delta^1_\alpha, \mu^2_\alpha.\delta^2_\alpha, \mu^1_\alpha.\delta^2_\alpha,
   \mu^2_\alpha.\delta^1_\alpha\} ,\\
    \max\{\mu^1_\alpha.\delta^1_\alpha ,\mu^2_\alpha.\delta^2_\alpha, \mu^1_\alpha.\delta^2_\alpha, \mu^2_\alpha.\delta^1_\alpha\}]$\\
    $[\mu\oslash\delta]_\alpha = [\min\{\mu^1_\alpha / \delta^1_\alpha, \mu^2_\alpha / \delta^2_\alpha, \mu^1_\alpha / \delta^2_\alpha,
   \mu^2_\alpha / \delta^1_\alpha\} ,\\
    \max\{\mu^1_\alpha / \delta^1_\alpha ,\mu^2_\alpha / \delta^2_\alpha, \mu^1_\alpha / \delta^2_\alpha, \mu^2_\alpha / \delta^1_\alpha\}]$\\
$[|\mu|]_\alpha=[\max \{0,\mu^1_\alpha,-\mu^2_\alpha\} , \max\{|\mu^1_\alpha| , |\mu^2_\alpha|\}]$
  $\label{7.2}$
\end{definition}
where$[\mu]_\alpha=[\mu^1_\alpha , \mu^2_\alpha],[\delta]_\alpha=[\delta^1_\alpha, \delta^2_\alpha]$

\begin{definition}
(\cite{M-M,K-S}) Let $\mu,\delta\in \mathcal{F}(\mathbb{R})$ and $[\mu]_\alpha=[\mu^1_\alpha, \mu^2_\alpha]$ and $[\delta]_\alpha=[\delta^1_\alpha ,
\delta^2_\alpha]$.\\
 The partial order $\preceq$ in $\mathcal{F}(\mathbb{R})$ is defined as below:\\
 $\mu\preceq\delta$ if and only if $\mu^1_\alpha \leq \delta^1_\alpha$ and $\mu^2_\alpha \leq \delta^2_\alpha$ for all $\alpha\in (0,1]$.\\
The equality of fuzzy numbers $\mu$ and $\delta$ is defined as below:\\
 $\mu = \delta \Leftrightarrow [\mu]_\alpha=[\delta]_\alpha$ for all $\alpha\in (0,1]$.
  $\label{7.3}$
\end{definition}

\subsection{Fuzzy Soft Real Number}
 \begin{definition}{\bf:}
 Let $\mathbb{R}$ be the set of all real numbers, and $E$ be the parameter set. The set of all fuzzy soft real numbers, say F.S real numbers, is denoted by
 $\mathbb{R}_E=\{f : E \rightarrow I^{\mathbb{R}}\}$, where for all $e\in E$, $f(e)$'s are fuzzy real numbers.
 $\label{def7.5}$
\end{definition}
%A typical fuzzy soft real number is shown in Fig. $\ref{Fig.4}$. For all parameters, the membership functions for this F.S real number is defined as a
%triangle membership function. But the shape of F.S real numbers are not limited to this graph.
% In other words, it is not necessary that the membership functions for all parameters follow the same rule.\\
  The real number $r\in \mathbb{R}$ can be seen as a F.S real number $\bar{r}_E$ if for all $e\in E$, $r(e)$ is defined by
characteristic function of $r$. So we have
 \begin{eqnarray*}
\bar{r} : E &\rightarrow& I^\mathbb{R}\\
 e &\mapsto& \bar{r}(e)=\chi_r
 \end{eqnarray*}
  where $\bar{r}(e)$ is defined by
$$\bar{r}_e (t)= \left\{\begin{array}{ll}
 1  & \mbox {if $t=r$}\\
 0 & \mbox{if $t\neq r$}\\
 \end{array}
 \right.$$
 for all $t\in \mathbb{R}$.
We call such a F.S real number, crisp F.S real number and denote it by $\bar{r}_E$. The F.S real number $f_E$ is called non-negative F.S real number, if
for all $e\in E$, $f_e(t)=0$ $\forall t<0$. The set of all non-negative F.S real numbers is denoted by $\mathbb{R}^*_E$.
 $\bar{0}_E$ and $\bar{1}_E$ are used to show the crisp F.S real numbers 0 and 1, respectively. So for
each $e\in E$ and for all $t\in \mathbb{R}$
 $$\bar{0}_e(t)= \left\{\begin{array}{ll}
 1  & \mbox {if $t=0$}\\
 0 & \mbox{if $t\neq 0$}\\
 \end{array}
 \right.$$
$$\bar{1}_e(t)= \left\{\begin{array}{ll}
 1  & \mbox {if $t=1$}\\
 0 & \mbox{if $t\neq 1$}\\
 \end{array}
 \right.$$
During this study, we will show the F.S real numbers by $\tilde r_E$, where $r\in \mathbb{R}$ and $\tilde r : E \rightarrow I^{\mathbb{R}}$.

%\begin{figure}
% \includegraphics[width=3.25 in]{Doc6}
% \caption{Triangle F.S real number in $\mathbb{R}_{(0,1]}$}\label{Fig.4}
%\end{figure}

\begin{remark}{\bf:}
In daily-life situations, saying like '' weather is pretty warm, although it is winter '' are used usually. So a mathematical concept is needed for
medelling these kinds of term.
\end{remark}

\begin{remark}{\bf:}
In application, F.S real numbers are used to display numbers which are approximately equal to a  real number or being approximately between two real
numbers whit respect to some parameters.
\end{remark}

\begin{definition}{\bf:}
$[\tilde r_E]_{e,\alpha}$ is called the $\alpha$-level set of F.S real number $\tilde r_E$ corresponding to the parameter $e\in E$ and defined as below
$$[\tilde r_E]_{e,\alpha}=\{t: \tilde r_e(t)\geq \alpha\}$$
So it can be considered as the $\alpha$-level set of fuzzy real number $\tilde r(e)$, i.e. $[\tilde r_E]_{e,\alpha}=[\tilde r(e)]_\alpha$.
 $\label{def7.5.1}$
\end{definition}

\begin{definition}{\bf:}
Let  $\tilde r_E$ and $\tilde r'_E$ be two F.S real numbers. We say that
\begin{enumerate}
    \item $\tilde r_E \tilde\preceq \tilde r'_E \Leftrightarrow \tilde r(e) \preceq \tilde r'(e)$, for all $e\in E$.
    \item  $\tilde r_E = \tilde r'_E \Leftrightarrow \tilde r(e) = \tilde r'(e)$, for all  $e\in E$.
\end{enumerate}
So if $$[\tilde r_E]_{e,\alpha}=[\tilde r (e)]_{\alpha}=[r^1_{e,\alpha} , r^2_{e,\alpha}]$$ and
  $$[\tilde r'_E]_{e,\alpha}=[\tilde r'(e)]_{\alpha}=[r'^1_{e,\alpha} , r'^2_{e,\alpha}]$$
   then we have
\begin{enumerate}
    \item $\tilde r_E \tilde\preceq \tilde r'_E \Leftrightarrow r^1_{e,\alpha} \leq r'^1_{e,\alpha}$ and $r^2_{e,\alpha} \leq r'^2_{e,\alpha}$
    $\forall \alpha\in (0,1]$ and $\forall e\in E$.
    \item  $\tilde r_E = \tilde r'_E \Leftrightarrow r^1_{e,\alpha} = r'^1_{e,\alpha}$ and $r^2_{e,\alpha} = r'^2_{e,\alpha}$,
    $\forall\alpha\in (0,1]$ and  $\forall e\in E$.
\end{enumerate}
 $\label{def7.5.2}$
\end{definition}

\subsection{Arithmetic Operations of Fuzzy Soft Real Number}
\begin{definition}{\bf:}
Let $\tilde r_E$ and $\tilde r'_E$ be two F.S real numbers. We define arithmetic operations $\tilde\oplus, \tilde\ominus, \tilde\otimes, \tilde\oslash$
 by applying Definition $\ref{7.1}$ as below:\\
$\tilde r_E \tilde\oplus \tilde r'_E$ is defined by map $\tilde r \tilde\oplus \tilde r': E \rightarrow I^\mathbb{R}$ such that for each $e\in E$ we have
 $(\tilde r \tilde\oplus \tilde r')(e)=\tilde r(e)
\oplus \tilde r'(e)$ where $\forall t\in\mathbb{R}$, $[\tilde r(e)
\oplus \tilde r'(e)](t)$ is defined by Definition $\ref{7.1}$,\\
$(\tilde r_E \tilde\ominus \tilde r'_E)(e)=\tilde r(e) \ominus \tilde r'(e)$, for each $e\in E$,\\
$(\tilde r_E \tilde\otimes \tilde r'_E)(e)=\tilde r(e) \otimes \tilde r'(e)$, for each $e\in E$,\\
$(\tilde r_E \tilde\oslash \tilde r'_E)(e)=\tilde r(e) \oslash \tilde r'(e)$, for each $e\in E$. $\label{def7.5.3}$
\end{definition}

 If $[\tilde r_E]_{e,\alpha}=[\tilde r(e)]_\alpha = [r^1_{e,\alpha} , r^2_{e,\alpha}]$, and $[\tilde r'_E]_{e,\alpha}=[\tilde r'(e)]_\alpha = [r'^1_{e,\alpha} , r'^2_{e,\alpha}]$
 be non-negative F.S real numbers,
then by applying Definition $\ref{7.2}$, we have\\
$[\tilde r_E\tilde \oplus \tilde r'_E]_{e,\alpha} = [\tilde r(e) \oplus \tilde r'(e)]_\alpha $, \\
$[\tilde r_E\tilde \ominus \tilde r'_E]_{e,\alpha} =[\tilde r(e) \ominus \tilde r'(e)]_\alpha $, \\
$[\tilde r_E\tilde \otimes \tilde r'_E]_{e,\alpha} =[\tilde r(e) \otimes \tilde r'(e)]_\alpha$, \\
$[\tilde r_E\tilde \oslash \tilde r'_E]_{e,\alpha} =[\tilde r(e) \oslash \tilde r'(e)]_\alpha$, \\
$[|\tilde r_E|]_{e,\alpha}=[|\tilde r(e)|]_\alpha$. \\
  $\forall e\in E$ and $\forall\alpha \in (0,1]$.

\begin{proposition}{\bf:}
The additive and the multiplicative identities in $\mathbb{R}_E$ are $\bar{0}_E$ and $\bar{1}_E$, respectively. $\label{pro7.1}$
 \end{proposition}
\textbf{Proof.}
 Let $r_E$ be a F.S real number and $\bar{0}_E$ be crisp F.S zero element. It is clear that
 $[\bar{0}_E]_{e,\alpha}=[\bar{0}(e)]_\alpha = \{0\}=[0,0]$.
 Suppose that $[\tilde r_E]_{e,\alpha}=[\tilde r(e)]_\alpha = [r^1_{e,\alpha} , r^2_{e,\alpha}]$.
 Then for
each $\alpha\in (0,1]$,
\begin{eqnarray*}
[r_E\tilde \oplus \bar{0}_E]_{e,\alpha} &=& [\tilde r(e) \oplus \bar{0}(e)]_\alpha \\
                                        &=& [r^1_{e,\alpha} +0 , r^2_{e,\alpha} +0]\\
                                        &=& [r^1_{e,\alpha} , r^2_{e,\alpha}] = [\tilde r(e)]_\alpha = [\tilde r_E]_{e,\alpha}
  \end{eqnarray*}
 Thus $\tilde r_E \tilde \oplus \bar{0}_E=\tilde r_E$.\\
 Similarly, it can be shown that $\forall e\in E$ and $\forall\alpha\in (0,1]$,
   $[\tilde r_E\tilde \otimes \bar{1}_E]_{e,\alpha} = [\tilde r_E]_{e,\alpha} \Leftrightarrow  \tilde r_E\tilde \otimes \bar{1}_E=\tilde r_E$.
$\Box$
%%%%%%%%%%%%%%%%%%%%%%%%%%%%%%%%%%%%%%%%%%%%%%%%%%%%%%%%%%%%%%%%%%%%%%%%%%%%%%%%%%%%%%%%%%%
\section{Fuzzy Soft Normed Spaces}
\subsection{Fuzzy Soft Norm}
 Due to introduce a norm for a fuzzy soft point, extension of some operations such as addition and scalar multiplication
are needed. In this section we firstly suggest a parametrization version of these well-known concepts, and then introduce the notion of fuzzy soft norm.
\begin{definition}{\bf:}
Let $f_E,g_E\in X_E$. Then we define $f_E \tilde\otimes g_E\in X_E\times X_E$ as the F.S multiplication by the map
\begin{eqnarray*}
 f \tilde\otimes g : E &\rightarrow& I^{X\times X}\\
                     e &\rightarrow& (f \tilde\otimes g)(e)=f(e) \otimes g(e)
 \end{eqnarray*}
 where
 \begin{eqnarray*}
(f \tilde\otimes g)(e)(x_1,x_2) &=& (f(e) \otimes g(e))(x_1,x_2)\\
                                &=& \min\{f_e(x_1) , g_e(x_2)\}
\end{eqnarray*}
 $\label{def7.1}$
\end{definition}

\begin{definition}{\bf:}
Let $X$ be a vector space over the field $F$ ($\mathbb{R}$ or $\mathbb{C}$).
 Let $h_1: X\times X \rightarrow X$ be the ordinary addition function on $X$, i.e.
$h_1(x_1,x_2)=x_1+x_2$ for all $x_1,x_2\in X$.\\
  Let $f_E,g_E\in X_E$. If F.S addition is denoted by $\tilde \oplus$, then by applying Definition $\ref{def7.0.1}$
   we have $h_1(f_E \times g_E)=f_E \tilde\oplus g_E$ where
\begin{eqnarray*}
(f\tilde \oplus g)(e)(x) &=& h_1(f \times g)(e)(x)\\
                         &=& \sup_{(x_1,x_2)\in h_1^{-1} (x)} [(f\times g)(e)](x_1,x_2)\\
                              &=& \sup_{(x_1,x_2)} [(f\times g)(e)](x_1,x_2)\\
                              &=& \sup_{(x_1,x_2)}\min\{f_e(x_1) , g_e(x_2)\}
\end{eqnarray*}
 where $x_1+x_2=x$.
 $\label{def7.2}$
\end{definition}

\begin{definition}{\bf:}
Let $X$ be a vector space over the field $F$ ($\mathbb{R}$ or $\mathbb{C}$).
 Let $h_2: F\times X \rightarrow X$ be the ordinary scalar multiplication function on $X$, i.e.
$h_2(t,x)=tx$ for $t\in F$ and $x\in X$.\\
 Let $f_E\in X_E$ and $r\in \mathbb{R}$ such that $\tilde r_E\in \mathbb{R}_E$. If F.S scalar
multiplication is denoted by $\tilde\otimes$, by applying Definition $\ref{def7.0.1}$ we have $h_2(r_E \times f_E)=r_E \tilde\otimes f_E$ where
\begin{eqnarray*}
(\tilde r\tilde \otimes f)(e)(z) &=& h_2(\tilde r\times f)(e)(z)\\
                                 &=& \sup_{(t,x)\in h^{-1}_2 (z)} [(\tilde r\times f)(e)](t,x)\\
                                      &=& \sup_{(t,x)} [(\tilde r\times f)(e)](t,x)\\
                                      &=& \sup_{(t,x)}\min\{\tilde r_e(t) , f_e(x)\}
\end{eqnarray*}
 where $tx=z$.
 $\label{def7.3}$
\end{definition}

Note that if $\tilde r_E$ be the crisp F.S real number $\bar{r}_E$, then the F.S scalar multiplication is defined by\\
 $(\bar{r} \tilde \otimes f)(e)(z)= \left\{\begin{array}{lll}
 f_e(\frac{z}{r}) & \mbox{if $r\neq 0$}\\
\sup_{x\in X} f_e(x) & \mbox{if $r=0,z=0$}\\
 0 & \mbox{if $r=0 , z\neq 0$}\\
 \end{array}
 \right.$\\

 We denote $\bar{r}_E \tilde \otimes f_E$ by $rf_E$.\\

 Next, let $f_E$ be the F.S point $\tilde x_E$, then\\
 $(\bar{r} \tilde \otimes \tilde{x})(e)(z)=
\left\{\begin{array}{ll}
 1 & \mbox{if $z=rx$}\\
 0 & \mbox{otherwise}\\
 \end{array}
 \right.$\\

So $\bar{r}_E \tilde \otimes \tilde{x}_E = {\overline{rx}}_E$, crisp F.S real number $rx$.

\begin{definition}{\bf:}
Let $X$ be a vector space over field $F$ ($\mathbb{R}$ or $\mathbb{C}$) and let $E$ be the set of parameters. We define fuzzy soft norm $||.||$, say
 F.S norm, over $X$ by map $||.|| : X_E \rightarrow \mathbb{R^*}_E$ which satisfies the below conditions
\begin{enumerate}
    \item $\tilde x_E = \bar{0}_E  \Leftrightarrow ||\tilde x_E||=\bar{0}_E$
    \item $||r\tilde x_E||= |\bar{r}_E| \tilde\otimes ||\tilde x_E||$
    \item $||\tilde x_E \tilde\oplus \tilde y_E|| \tilde \preceq ||\tilde x_E|| \tilde\oplus ||\tilde y_E||$
\end{enumerate}
$\label{7.6}$
\end{definition}
We denote the fuzzy soft normed space $X$, say F.S normed space $X$, by $(X,E,||.||)$.

 Note that F.S norm of F.S point $\tilde x_E \in X_E$
 is denoted by map $||\tilde x_E|| : E \rightarrow I^{\mathbb{R^*}}$ which is a non-negative F.S real number.
 It means that for each $e\in E$, $||\tilde x_E||(e) :
\mathbb{R} \rightarrow I$ where $||\tilde x_E||(e)(t)=0$ for all $t<0$ (see Definition $\ref{7.0}$).

 We can also define $\alpha$-level sets of F.S real number $||\tilde x_E||$ as below
  $$[||\tilde x_E||]_{e,\alpha} = [||\tilde x_E||(e)]_\alpha = [||\tilde x_E||^1_{e,\alpha} ,||\tilde x_E||^2_{e,\alpha}]$$
where $||\tilde x_E||^1_{e,\alpha} ,||\tilde x_E||^2_{e,\alpha}$ are real numbers and $[||\tilde x_E||(e)]_\alpha$, is a closed interval in real line (see
Definition $\ref{7.0}$).

\begin{lemma}{\bf:}
Let $r,x\in\mathbb{R}$. Suppose $\bar x_E$ and $\bar{r}_E$ be crisp F.S real numbers corresponding to $x$ and $r$, respectively.
 Then
\begin{enumerate}
    \item $|\bar {r}_E| = \overline{|r|}_E$
    \item ${\overline{rx}}_E=\bar{r}_E \tilde\otimes \bar{x}_E$
\end{enumerate}
   $\label{lem7.1}$
 \end{lemma}
\textbf{Proof.}
\begin{enumerate}
    \item Let $r\in \mathbb{R}$. Definition $\ref{def7.5}$ implies $[\bar{r}_E]_{e,\alpha}=[r , r]$
    and $[\overline{|r|}_E]_{e,\alpha}=[|r| , |r|]$. By applying Definition $\ref{def7.5.2}$, we have
\begin{eqnarray*}
  [|\bar{r}_E|]_{e,\alpha} &=& [|\bar {r}(e)|]_\alpha \\
                           &=& [\max \{0,r, -r\} , \max \{|r| , |r|\}]=[|r| , |r|]
\end{eqnarray*}
 So $|\bar {r}_E| = \overline{|r|}_E$.
    \item    It is clear by Definition $\ref{def7.3}$.
\end{enumerate}
$\Box$

\begin{theorem}{\bf:}
Let $(X,E,||.||)$ be a F.S normed space. Then for each $e\in E$,
\begin{enumerate}
    \item $||.||_e$ is a fuzzy norm, introduced in \cite{M-M,X-Z,x-z}, such that
      $|| . ||_e : X \rightarrow \mathcal{F}(\mathbb{R^*})$ is defined as
    $||x||_e =||\overline{x}_E||(e)$ where $\mathcal{F}(\mathbb{R^*})$ denotes the set of all non-negative fuzzy real numbers.
    \item   $(X,||.||^i_{e,\alpha})$ is
     normed space where $i=1,2$, $e\in E$ and $\alpha\in (0,1]$ as below\\
        $|| . ||^i_{e,\alpha} : X \rightarrow \mathbb{R^*}$ is defined by $\alpha$-level set of fuzzy norm $||.||_e$
    as $||x||^i_{e,\alpha} ={||x||_e}^i_\alpha$ where
      $$[||\bar x_E||]_{e,\alpha}=[||\bar x_E||(e)]_\alpha=[||x||_e]_\alpha =[{||x||_e}^1_\alpha , {||x||_e}^2_\alpha]$$
    \end{enumerate}
$\label{th7.1}$
\end{theorem}
\textbf{Proof.} Let $(X,E,||.||)$ be a F.S normed space.
\begin{enumerate}
    \item Take $e\in E$.
\begin{description}
    \item[(1)] Let $x=0$ and $\bar{0}$ denotes crisp fuzzy number zero. Then by applying Definitions $\ref{7.6}$, we have
    $\bar x_E = \bar{0}_E \Leftrightarrow ||\bar x_E||=\bar{0}_E
    \Leftrightarrow \forall e\in E, \forall t\in \mathbb{R}: ||\bar x_E||(e)=\bar{0}\Leftrightarrow ||x||_e =\bar{0}$.
    \item[(2)] Lemma $\ref{lem7.1}$ implies
    \begin{eqnarray*}
  [||rx||_e]_\alpha &=& \{ t: ||rx||_e(t) \geq \alpha\}\\
           &=& \{ t: ||{\overline{r x}}_E||)(e)(t)\geq \alpha \}\\
           &=& \{ t: ||r\bar x_E||)(e)(t)\geq \alpha \}\\
           &=& \{ t: [|\bar r_E| \tilde\otimes ||\bar x_E||](e)(t)\geq \alpha \}\\
           &=& \{ t: [\overline {|r|}_E \tilde\otimes ||\bar x_E||](e)(t)\geq \alpha \}\\
           &=& \{ t: [|r| ||\bar x_E||](e)(t)\geq \alpha \}\\
           &=& \{ t: ||\bar x_E||(e)(\frac{t}{|r|})\geq \alpha \}\\
           &=& \{ t: ||x||_e (\frac{t}{|r|})\geq \alpha \}
    \end{eqnarray*}
    On the other hand
 \begin{eqnarray*}
  [|r| ||x||_e]_\alpha &=& \{ t: [|r| ||x||_e](t) \geq \alpha\}\\
                       &=& \{ t: ||x||_e (\frac{t}{|r|})\geq \alpha \}
    \end{eqnarray*}
    \item[(3)] For any pair $\bar x_E,\bar y_E\in \tilde X$,
    $$||\bar x_E \tilde\oplus \bar y_E||\tilde \preceq ||\bar x_E|| \tilde\oplus ||\bar y_E||$$
    Then for each $e\in E$ and $\alpha\in (0,1]$,\\
    $$[||\bar x_E\tilde\oplus \bar y_E||(e)]^i_\alpha \leq [||\bar x_E||(e)]^i_\alpha + [||\bar y_E||(e)]^i_\alpha$$
when $i=1,2$. So
    $${||x+ y||_e}^i_\alpha \leq {||x||_e}^i_\alpha + {||y||_e}^i_\alpha$$
    where $x+y=$ supp $(\bar x_E\tilde\oplus \bar y_E)$.
         This implies that
     $||x+y||_e\preceq ||x||_e \oplus ||y||_e$.
     Thus $||.||_e$'s are fuzzy norm on $X$ for all $e\in E$.
\end{description}
    \item It is similar to (1).
     \end{enumerate}

\begin{theorem}{\bf:}
Let $(X,||.||)$ be a normed space. Then
\begin{enumerate}
    \item $||.||^{\mathcal{F}}$ is a fuzzy norm on $X$, where $||.||^{\mathcal{F}}: X\rightarrow \mathcal{F}(\mathbb{R^*})$ and for each $x\in X$,
     $||x||^{\mathcal{F}}=\chi_{||x||}$ where $\chi$ is a characteristic function.  So for all $t\in \mathbb{R^*}$
          $$ ||x||^{\mathcal{F}}(t) = \left\{\begin{array}{ll}
 1  & \mbox {if $t=||x||$}\\
 0 & \mbox{otherwise}\\
 \end{array}
 \right.$$
    \item $||.||^E$ is a F.S norm on $X$ with respect to the parameter set $E$. $||.||^E$ is defined by the mapping
    $||.||^E:  X_E \rightarrow \mathbb{R^*}_E$ such as for all $\tilde x_E\tilde\in X_E$,and $\forall e\in E$, $\forall t\in \mathbb{R^*}$
    $$(||\tilde x_E||^E)_e(t)= \left\{\begin{array}{ll}
 1  & \mbox {if $t=||x||$}\\
 0 & \mbox{otherwise}\\
 \end{array}
 \right.$$
    where  $x=$support $\tilde x(e)$ for all $e\in E$.
    \end{enumerate}
$\label{th7.2}$
\end{theorem}

\textbf{Proof.} Part (1) is clear. We only prove (2).
    \begin{description}
    \item[(1)] Let $\tilde x_E = \bar{0}_E$. Suppose $x=$support $\tilde x(e)$  for each $e\in E$, then for $e\in E$ we have
     \begin{eqnarray*}
      \tilde x_E = \bar{0}_E &\Leftrightarrow & \tilde x(e)=\bar{0}\\
                             &\Leftrightarrow & x=0\\
                             & \Leftrightarrow & ||x||=0\\
                             & \Leftrightarrow & (||\tilde x_E||^E)_e(0)=1, (||\tilde x_E||^E)_e(t)=0; t\neq 0\\
                             & \Leftrightarrow & ||\tilde x_E||^E=\bar{0}_E
      \end{eqnarray*}
    \item[(2)] It is clear that $supp (\bar{r}_E \tilde\otimes \tilde x_E)(e)=rx$, then
 \begin{eqnarray*}
      [||\bar{r}_E \tilde\otimes \tilde x_E||^E]_{e,\alpha} &=& [||\bar{r}_E \tilde\otimes \tilde x_E||^E (e)]_\alpha\\
                                                            &=& \{t: (||\bar{r}_E \tilde\otimes \tilde x_E||^E (e))(t)\geq \alpha\}\\
                                                            &=& \{t: t= || supp (\bar{r}_E \tilde\otimes \tilde x_E)(e) || \}\\
                                                            &=& \{t: t=||rx||\}=||rx||=|r| ||x||
      \end{eqnarray*}
      and
      \begin{eqnarray*}
      [\overline{|r|}_E \tilde\otimes ||\tilde x_E||^E]_{e,\alpha} &=& [(\overline{|r|}_E)(e) \otimes (||\tilde x_E||^E )(e)]_\alpha\\
                                                            &=& [ |r| ||x|| , |r| ||x|| ] = |r| ||x||\\
                                                                  \end{eqnarray*}
    \item[(3)] Let $\tilde x_E$, and $\tilde y_E\in \tilde X$. Then
    \begin{eqnarray*}
      [||\tilde x_E \tilde\oplus \tilde y_E||^E]_{e,\alpha} &=& [||\tilde x_E \tilde\oplus \tilde y_E||^E (e)]_\alpha\\
                                                            &=& \{t: (||\tilde x_E \tilde\oplus \tilde y_E||^E (e))(t)\geq \alpha\}\\
                                                            &=& \{t: t= || supp (\tilde x_E \tilde\oplus \tilde y_E)(e) || \}\\
                                                            &=& \{t: t=||x+y||\}=||x+y||
      \end{eqnarray*}
      and
\begin{eqnarray*}
      [||\tilde x_E||^E \tilde\oplus ||\tilde y_E||^E]_{e,\alpha} &=& [||\tilde x_E||^E(e) \oplus ||\tilde y_E||^E (e)]_\alpha\\
                                                            &=& [||x||+||y|| , ||x||+||y||]\\
                                                            &=& ||x||+||y||
      \end{eqnarray*}
      Hence $||\tilde x_E \tilde\oplus \tilde y_E||^E \tilde \preceq ||\tilde x_E||^E \tilde\oplus ||\tilde y_E||^E$.
\end{description}

\subsection{Fuzzy Soft Topology Generated by F.S Norm}
 If $(X,E,||.||)$ is a F.S normed space, then the F.S norm $||.||$ induces a F.S topology over $X$ as below.

\begin{definition}{\bf:}
Let $(X,E,||.||)$ be a F.S normed space. The F.S topological space $(X,E,w'(w(\tau_{||.||^i_{e,\alpha}})))$ is called F.S normed topology over $X$
generated by F.S norm $||.||$.
 $\label{7.11}$
\end{definition}

\begin{theorem}{\bf:}
Let $(X,E,||.||)$ be a F.S normed space. The topological space $(X,E,w'(w(\tau_{||.||^i_{e,\alpha}})))$ is a F.S Hausdorff space.
  $\label{th7.10}$
\end{theorem}

\textbf{Proof.} Let $\tilde  x_E$ and $\tilde y_E$ be two F.S disjoint points in F.S normed space $(X,E,||.||)$, where for all $e\in E$, $\tilde
x(e)=x^{\lambda_e}$ and $\tilde y(e)=y^{\gamma_e}$. Hence $x$ and $y$ are disjoint points in normed space $(X,||.||^i_{e,\alpha})$. Since every normed
space is Hausdorff, then there exist open sets $U$ and $V$ in $\tau_{||.||^i_{e,\alpha}}$ such that $U \cap V = \emptyset$. Consider F.S open sets $U_E$
and $V_E$ in $w'(w(\tau_{||.||^i_{e,\alpha}}))$
 such that for all $e\in E$, $U(e)=\chi_{U}$ and $V(e)=\chi_{V}$. It is clear
that $\tilde x_E\tilde \in U_E$ and $\tilde y_E\tilde \in V_E$ and moreover $U_E \tilde \wedge V_E = \Phi$. $\Box$

\subsubsection{Convergency in Fuzzy Soft Normed Spaces}
 Here we suggest the concept of ''approach'' and consequently ''limit'' in the fuzzy soft set theory.
  In point-set topology, $\lim_{x\rightarrow c} f(x) = L$ means that when $x$ approaches $c$ sufficiently, $f(x)$ becomes
  arbitrarily close to $L$.\\
  By extending this idea in fuzzy soft set theory, where a typical fuzzy soft set is in fact a map whose image is a function from $X$ into
  $[0,1]$, the concept of ''approximately near with respect to some parameters'' can be presented as the following.\\
  Let $\tilde x_E$ and $\tilde y_E$ be F.S points in $X$.
  ''$\tilde x_E$ approaches $\tilde y_E$'' means that,
 for all $e\in E$ the real function $\tilde x(e)$  becomes close to the real function $\tilde y(e)$ whenever
  $x$ approaches $y$ in normed space $(X, ||.||^i_{e,\alpha})$. In other words,
 behavior of objects $x$ and $y$
on the basis of some parameters are similar.

\begin{definition}{\bf:}
Let $(X,E,||.||)$ be a F.S normed space. Let $\tilde x_E$ and $\tilde y_E$ be two F.S points in $X$ such that for each $e\in E$
   $$\tilde x_e(z)= \left\{\begin{array}{ll}
 \lambda_e & \mbox{if $z=x$}\\
 0 & \mbox{if $z\neq x$}\\
 \end{array}
 \right.$$
  $$\tilde y_e(z)= \left\{\begin{array}{ll}
 \gamma_e & \mbox{if $z=y$}\\
 0 & \mbox{if $z\neq y$}\\
 \end{array}
 \right.$$
 where $\lambda_e , \gamma_e \in (0,1]$. We say that $\tilde x_E$ approaches  $\tilde y_E$,
  denoting by
   $$\tilde x_E \rightarrow \tilde y_E$$
    whenever for all $e\in E$ and $\forall \varepsilon >0$
 there exists $\delta >0$ such that
 $$||x-y||^i_{e,\alpha} < \delta \Rightarrow |\lambda_e - \gamma_e| < \varepsilon$$
 for $i=1,2$ and $\forall\alpha\in (0,1]$.
  $\label{7.7.1}$
\end{definition}

\begin{definition}{\bf:}
Let $(X,E,||.||)$ be a F.S normed space. Let $\{x_n\}$ be a sequence in $X$ where $n=1,2,3,\ldots$ .
 For each $n$, let $(\tilde x_E)_n$ be a F.S point in $X$ such that
  $\forall e\in E$ and $\forall z\in X$
 $$(\tilde x_n)_e(z)= \left\{\begin{array}{ll}
 \lambda_{e,n} & \mbox{if $z=x_n$}\\
 0 & \mbox{if $z\neq x_n$}\\
 \end{array}
 \right.$$
 where  $\lambda_{e,n} \in (0,1]$. Then $\{ (\tilde x_E)_n \}$ includes F.S points $(\tilde x_E)_n $ is called
  a F.S sequence in $X$.
 $\label{7.7.2}$
\end{definition}

\begin{definition}{\bf:}
Let $(X,E,||.||)$ be a F.S normed space. Let $\{ (\tilde x_E)_n \}$ be a F.S sequence in $X$ mentioned in Definition $\ref{7.7.2}$.
 Suppose that $\tilde x_E$ be a F.S point in $X$, such that $\forall e\in E$ and $\forall z\in X$
 $$\tilde x_e(z)= \left\{\begin{array}{ll}
 \gamma_e & \mbox{if $z=x$}\\
 0 & \mbox{if $z\neq x$}\\
 \end{array}
 \right.$$
 where $\gamma_e\in (0,1]$. By using Definition $\ref{7.7.1}$, we say that the sequence $\{ (\tilde x_E)_n \}$
  converges to $\tilde x_E$ denoting by
 $$(\tilde x_E)_n \rightarrow \tilde x_E$$
    if and only if
   $\forall \varepsilon >0$
 there exist positive integer $N$ and $\delta >0$ such that for all $e\in E$,
 $$\forall n\geq N; ||x_n-x||^i_{e,\alpha} < \delta \Rightarrow |\lambda_{e,n} - \gamma_e| < \varepsilon$$
 or  $\lim_{n\rightarrow\infty} \lambda_{e,n} = \gamma_e$ whenever in normed spaces $(X,||.||_{e,\alpha}^i)$
  we have $\lim_{n\rightarrow\infty} x_n = x$.
  We call the sequence $\{ (\tilde x_E)_n \}$ a F.S convergent sequence in $X$.
 $\label{7.7}$
\end{definition}

\begin{theorem}{\bf:}
Let $(X,E,||.||)$ be a F.S normed space. If F.S sequence $\{ (\tilde x_E)_n \}$ converges to $\tilde x_E$ and $\tilde y_E$, then $\tilde x_E=\tilde y_E$.
$\label{th7.5}$
\end{theorem}

\textbf{Proof.} Let $\{ (\tilde x_E)_n \}$ be a F.S sequence in $X$ mentioned in Definition $\ref{7.7.2}$. Let $(\tilde x_E)_n \rightarrow \tilde x_E$ and
$(\tilde x_E)_n \rightarrow \tilde y_E$ where $\tilde x_E$ and $\tilde y_E$ are two differen F.S points in $X$ such that $\forall e\in E$, are defined as
below
$$\tilde x_e(z)= \left\{\begin{array}{ll}
 \xi_e  & \mbox {if $z=x$}\\
 0 & \mbox{otherwise}\\
 \end{array}
 \right.$$
$$\tilde y_e(z)= \left\{\begin{array}{ll}
 \gamma_e  & \mbox {if $z=y$}\\
 0 & \mbox{otherwise}\\
 \end{array}
 \right.$$
 where $z\in X$ and $\xi_e,\gamma_e \in (0,1]$.\\
   Since limit is unique in every normed (or metric) spaces, Definition $\ref{7.7}$ implies that
   $x=y$ and $\xi_e=\gamma_e $, $\forall e\in E$. This means that $\tilde x_E=\tilde y_E$. $\Box$

\begin{definition}{\bf:}
Let $(X,E,||.||)$ be a F.S normed space. Let $\{(\tilde x_E)_n \}$ be a F.S sequence in $X$ mentioned in Definition $\ref{7.7.2}$.
   We say that the F.S sequence $\{ (\tilde x_E)_n \}$ is a fuzzy soft Cauchy sequence, say F.S Cauchy sequence,  in $X$ if F.S points $(\tilde x_E)_n$ become
    arbitrarily close to each other as the sequence progress.\\
     So we say $\{ (\tilde x_E)_n \}$ is a F.S Cauchy sequence in $X$
   if and only if
   $\forall \varepsilon >0$ there exists a positive integer $N$ such that $\forall n,m \geq N$
    $$|\lambda_{e,n}- \lambda_{e,m}| \leq \varepsilon$$
     while $$|| x_n - x_m||_{e,\alpha}^i \leq \varepsilon$$
      $\label{7.8}$
\end{definition}

\begin{definition}{\bf:}
Let $(X,E,||.||)$ be a F.S normed space and $\{(\tilde x_E)_n \}$ be a F.S sequence in $X$. Let $\{n_k\}$ be a sequence of positive integers such that
$n_1<n_2<\ldots$.
   Then the F.S sequence $\{(\tilde x_{E})_{n_k} \}$ is called a subsequence of $\{(\tilde x_E)_n\}$.
 $\label{7.9}$
\end{definition}

\begin{theorem}{\bf:}
Let $(X,E,||.||)$ be a F.S normed space and $\{ (\tilde x_E)_n \}$ be a F.S sequence in $X$. The sequence $\{ (\tilde x_E)_n\}$ is converges to $\tilde
y_E$ in $X$ if and only if every subsequence of it is convergent to $\tilde y_E$. $\label{th7.3}$
\end{theorem}

\textbf{Proof.} Let $\{ (\tilde x_E)_n \}$ be a F.S sequence mentioned in Definition $\ref{7.7.2}$, and $\tilde y_E$ be a F.S point in $X$ such that
$\forall e\in E$ and $\forall z\in X$,
 $$\tilde y_e(z)= \left\{\begin{array}{ll}
 \gamma_e  & \mbox {if $z=y$}\\
 0 & \mbox{otherwise}\\
 \end{array}
 \right.$$
 where $\gamma_e \in (0,1]$.
\begin{description}
    \item[$\Rightarrow$] Let the sequence $\{(\tilde x_E)_n \}$ converges to $\tilde y_E$.
      Definition $\ref{7.7}$ implies that for each $e\in E$, sequence $\{\lambda_{e,n}\}$ of real numbers converges to $\gamma_e$
      while sequence $\{x_n\}$ converges to $y$
      in normed spaces $(X,||.||_{e,\alpha}^i)$.
       So for every subsequence $\{n_k\}$ of $\{n\}$, we have
       $$\lim_{n_k\rightarrow\infty} \lambda_{e,n_k}=\gamma_e$$
         whenever $$\lim_{n_k\rightarrow\infty} x_{n_k} = y$$
          This implies that
     $(\tilde x_{E})_{n_k} \rightarrow \tilde y_E$.
    \item[$\Leftarrow$] Let $\{n_k\}$ be an arbitrary subsequence of $\{n\}$. So the subsequence
     $\{(\tilde x_{E})_{n_k}\}$ of $\{(\tilde x_E)_n \}$ is convergent to $\tilde y_E$.
 Definition $\ref{7.7}$ implies $\lambda_{e,n_k} \rightarrow \gamma_e$ whenever $x_{n_k} \rightarrow y$ in  normed spaces $(X,||.||_{e,\alpha}^i)$.
     Since $\{n_k\}$ is an arbitrary subsequence of $\{n\}$, then
        $\lambda_{e,n} \rightarrow \gamma_e$ while $x_{n} \rightarrow y$.
 Hence $(\tilde x_E)_n \rightarrow \tilde y_E$. This complete the proof. \\
 $\Box$
\end{description}

\begin{theorem}{\bf:}
Let $(X,E,||.||)$ be a F.S normed space. Then every F.S convergent sequence is a F.S Cauchy sequence. $\label{th7.4}$
\end{theorem}

\textbf{Proof.} It is implied from Definitions $\ref{7.7}$ and $\ref{7.8}$.\\
$\Box$

\begin{definition}{\bf:}
 Let $(X,E,||.||)$ be a F.S normed space. If every F.S Cauchy sequence in $X$ be a F.S convergent sequence, the F.S normed space $X$ is called a F.S Banach
 space.
 $\label{7.16}$
\end{definition}

\subsubsection{Continuity in Fuzzy Soft Normed Spaces}
\begin{definition}{\bf:}
Let $X$ and $Y$ be universal sets and $E$ and $E'$ are the parameter sets. Let $(X,E,||.||)$ and $(Y,E',||.||)$ be F.S normed spaces.
 Then the map $$T: (X,E,||.||) \rightarrow (Y,E',||.||)$$ is called a fuzzy soft operator, say F.S operator in brief.
 $\label{7.10}$
\end{definition}

\begin{definition}{\bf:}
Let  $(X,E,||.||)$ and $(Y,E',||.||)$ be F.S normed spaces where $E$ and $E'$ are the parameter sets and $X$ and $Y$ are universal sets.
 Then $T: (X,E,||.||) \rightarrow (Y,E',||.||)$ is continuous at F.S point $\tilde x_E \tilde\in X_E$, called F.S continuous, if
    for every $\tilde z_E \tilde \in X_E$, such that
  $\tilde z_E \rightarrow \tilde x_E$ we have  $T \tilde z_E \rightarrow  T \tilde x_E$.
 $\label{7.13}$
\end{definition}
We say that $T$ is F.S continuous over $X$ whenever $T$ is continuous at every F.S point of $X$.

\begin{lemma}{\bf:}
Let $(X,E,||.||)$ and $(Y,E',||.||)$ be F.S normed spaces. Let $T: (X,||.||^i_{e,\alpha}) \rightarrow (Y,||.||^i_{e',\alpha})$ be an operator from $X$ into
$Y$. Then if there is a bijection between $E$ and $E'$, then for every F.S point $\tilde x_E$ in $X$ $T\tilde x_E={\widetilde{Tx}}_{E'}$.
  $\label{lem5.1}$
\end{lemma}
\textbf{Proof.} Definition $\ref{def7.0.1}$ implies that
 $$(T\tilde x)_{e'}(y) = \left\{\begin{array}{ll}
 \tilde x_e(x)  & \mbox {if $y=Tx$}\\
 0 & \mbox{otherwise}\\
 \end{array}
 \right.$$
 This means that $T\tilde x_E$ is a F.S point in $Y$ corresponding to $Tx\in Y$. This complete the proof.\\
$\Box$

\begin{theorem}{\bf:}
Let $(X,E,||.||)$ and $(Y,E',||.||)$ be two F.S normed spaces. Let $T: (X,||.||^i_{e,\alpha}) \rightarrow (Y,||.||^i_{e',\alpha})$ be an operator from $X$
into $Y$. Let $\tilde x_E$ and $\tilde z_E$ are used to denote F.S points in $X$ corresponding to $x,z\in X$, respectively. Then $T$ is continuous at
$\tilde x_E\tilde\in X_E$ if and only if $T$ be continuous at $x\in X$
 and moreover $\lim_{z\rightarrow x} |\tilde x_e(x)-\tilde z_e(x)|=0$ for all $z\in X$.
  $\label{th7.7.1}$
\end{theorem}
\textbf{Proof.} It is clear by Lemma $\ref{lem5.1}$ and Definitions $\ref{7.7}$ and $\ref{7.13}$.\\
$\Box$

\begin{theorem}{\bf:}
Let $(X,E,||.||)$ and $(Y,E',||.||)$ be two F.S normed spaces and $T: (X,||.||^i_{e,\alpha}) \rightarrow (Y,||.||^i_{e',\alpha})$ be an operator from $X$
into $Y$. Let $g_{E'}$ be a F.S open subset of $\tilde Y$ with respect to norm topology $w'(w(\tau_{||.||^i_{e',\alpha}}))$ induced by F.S norm $||.||$
where $i=1,2$ and $\alpha\in (0,1]$. Then $T$ is F.S continuous over $X$ if and only if $T^{-1}(g_{E'})$ is F.S open subset of $\tilde X$ with respect to
norm topology $w'(w(\tau_{||.||^i_{e,\alpha}})))$.
   $\label{th7.7}$
\end{theorem}
\textbf{Proof.} It is clear by applying Theorem $\ref{th7.7.1}$. $\Box$

\begin{definition}{\bf:}
 $T: (X,E,||.||) \rightarrow (Y,E',||.||)$ is called F.S sequentially continuous if whenever the sequence $\{(\tilde x_E)_n\}$ of F.S points in $X$
  converges to $\tilde x_E$, the F.S sequence $\{T(\tilde x_E)_n\}$ converges to F.S point $T\tilde x_E$.
 $\label{7.14}$
\end{definition}

\begin{theorem}{\bf:}
Let $(X,E,||.||)$ and $(Y,E',||.||)$ be two F.S normed spaces and $T: (X,||.||^i_{e,\alpha}) \rightarrow (Y,||.||^i_{e',\alpha})$ be an operator from $X$
into $Y$. Then $T$ is F.S sequentially continuous if and only if $T$ is F.S continuous over $X$.
  $\label{th7.8}$
\end{theorem}

\textbf{Proof.} It is clear by applying Theorem $\ref{th7.7.1}$. $\Box$
%%%%%%%%%%%%%%%%%%%%%%%%%%%%%%%%%%%%%%%%%%%%%%%%%%%%%%%%%%%%%

   \section{Fuzzy Soft Fixed-Point Theorem}
   In this section we give a parameterized extension of fixed-point theorem in fuzzy soft set normed spaces.

 \begin{definition}{\bf:}
Let $X$ and $E$ are sets of objects and parameters, respectively. Let $T: (X,||.||^i_{e,\alpha}) \rightarrow (X,||.||^i_{e,\alpha})$ be an operator. The
F.S point $\tilde x_E$ is called
 a F.S fixed point of $T$ if and only if $T\tilde x_E= \tilde x_E$.
  $\label{th7.15}$
\end{definition}

 \begin{theorem}{\bf:}
Let $(X,E,||.||)$ be a F.S Banach spaces. Let $T: (X,||.||^i_{e,\alpha}) \rightarrow (X,||.||^i_{e,\alpha})$ be an contraction operator, i.e., for some
real number like $k$ such that $0<k<1$ we have
 $$|| T x - T y||^i_{e,\alpha} \leq k || x - y||^i_{e,\alpha}$$
 where $i=1,2$ , $e\in E$ ,$\alpha\in (0,1]$ and $x,y\in X$.
 Then $T$ is a F.S continuous map on $X$ and moreover
  $$|| T \tilde x_E \tilde\ominus T \tilde y_E || \tilde\preceq k || \tilde x_E \tilde\ominus \tilde y_E||$$
  $\label{th7.9}$
\end{theorem}
\textbf{Proof.} It is clear by applying Theorem $\ref{th7.7.1}$ and Definition $\ref{def7.5.2}$. $\Box$

 \begin{theorem}{\bf:}
Let $(X,E,||.||)$  be a F.S Banach spaces. Let $T: (X,||.||^i_{e,\alpha}) \rightarrow (X,||.||^i_{e,\alpha})$ be an contraction operator, i.e., for some
real number like $k$ such that $0<k<1$
 $$|| T x - T y||^i_{e,\alpha} \leq k || x - y||^i_{e,\alpha}$$
 where $i=1,2$ and $e\in E$ and $\alpha\in (0,1]$.
 Then there exists a unique F.S point in $X$ like $\tilde x_E$ such that $T\tilde x_E= \tilde x_E$. Moreover for any
  $\tilde z_E$ in $X$, $T(T(\ldots (T\tilde z_E))\ldots) \rightarrow \tilde x_E$.
  $\label{th7.10}$
\end{theorem}

\textbf{Proof.} Let $X$ and $E$ be the sets of objects and parameters, respectively. Let $\tilde z_E$ be a F.S point in $X$ such that
  $$(\tilde z)_{e}(x) = \left\{\begin{array}{ll}
 \xi_e  & \mbox {if $x=z$}\\
 0 & \mbox{otherwise}\\
 \end{array}
 \right.$$
 where for all $e\in E$, $\xi_e\in (0,1]$. Define the F.S sequence
$\{(\tilde x_E)_n\}$ in $X$ as below:\\
 $T \tilde z_E = (\tilde x_E)_1$\\
 $T (\tilde x_E)_1 = (\tilde x_E)_2$\\
 $ T (\tilde x_E)_2 = (\tilde x_E)_3$\\
   $\vdots$\\
   $T (\tilde x_E)_{n-1} = (\tilde x_E)_n$\\
    $\vdots$\\
 where $\tilde x_n: E \rightarrow I^X$
 such that
 $$(\tilde x_n)_{e}(x) = \left\{\begin{array}{ll}
 \lambda_{e,n}  & \mbox {if $x=x_n$}\\
 0 & \mbox{otherwise}\\
 \end{array}
 \right.$$
   Since $\tilde z_E$ and $(\tilde x_E)_n$'s are F.S points in $X$, for all $e\in E$, we can define the sequence $\{x_n\}$ in $X$ as below:\\
 $x_1 = supp \tilde x_1 (e)=Tz$\\
 $x_2 = supp \tilde x_2 (e)=Tx_1=T(Tz)=T^2 z$\\
    $\vdots$\\
   $x_n = supp \tilde x_n (e)=Tx_{n-1}=T(T(\ldots (Tz))\ldots)=T^n z$\\
    $\vdots$\\
    where for every $i=1,2,\ldots$, $supp \tilde x_i (e) =[\tilde x_i (e)]^{-1}(0,1]$  indicates the support of fuzzy set $\tilde x(e)$ and $T^n$
     shows the composition of $T$ with itself $n$
    times.\\
    On the other hand, Lemma $\ref{lem5.1}$ implies that the membership functions of $(\tilde x_E)_n$'s are defined as below:
 $$(\tilde x_1)_{e}(x) = (T\tilde z)_e(x) = \left\{\begin{array}{ll}
 \xi_e  & \mbox {if $x=x_1=z$}\\
 0 & \mbox{otherwise}\\
 \end{array}
 \right.$$
$$(\tilde x_2)_{e}(x) = (T\tilde x_1)_e(x) = \left\{\begin{array}{ll}
 (\tilde x_1)_{e}(x)=\xi_e  & \mbox {if $x=x_2=x_1=z$}\\
 0 & \mbox{otherwise}\\
 \end{array}
 \right.$$
 $\vdots$
$$(\tilde x_n)_{e}(x) = (T\tilde x_{n-1})_e(x) = \left\{\begin{array}{ll}
 \xi_e  & \mbox {if $x=x_n=\ldots=x_1=z$}\\
 0 & \mbox{otherwise}\\
 \end{array}
 \right.$$
$\vdots$\\
So for all $e\in E$ we have $\lambda_{e,n} =\xi_e $, $\forall n\in \mathbb{N}$.\\
     Since $T$ is a contraction map, we can show that
 $$|| x_{n+1}- x_n||^i_{e,\alpha} \leq k^n || x_1 - z||^i_{e,\alpha}$$
 and moreover, for any $\varepsilon > 0$ since $0< k <1$ we can fined a large number
 $N\in \mathbb{N}$ such that $k^N < \frac{\epsilon (1-k)}{ || x_1 - z||^i_{e,\alpha}}$. So for all $m,n\in \mathbb{N}$ such that $m >
 n\geq N$ we have
 \begin{eqnarray*}
 || x_m - x_n||^i_{e,\alpha} & \leq & ||x_m - x_{m-1}||^i_{e,\alpha} + ||x_{m-1}-x_{m-2}||^i_{e,\alpha}\\
                             &      & \mbox{} + ... +||x_{n+1}-x_n||^i_{e,\alpha}\\
                             & \leq & k^{m-1} ||x_1-z||^i_{e,\alpha} +k^{m-2} ||x_1-z||^i_{e,\alpha}\\
                             &      & \mbox{}  +...+k^n ||x_1-z||^i_{e,\alpha}\\
                             & = & ||x_1-z||^i_{e,\alpha} k^n \sum_{i=0}^{m-n-1} k^i\\
                             & \leq & ||x_1-z||^i_{e,\alpha} k^n \sum_{i=0}^{\infty} k^i \\
                             & = & ||x_1-z||^i_{e,\alpha} k^n \frac{1}{1-k}\\
                             & \leq & \frac{\epsilon(1-k)}{||x_1-z||^i_{e,\alpha}}.\frac{||x_1-z||^i_{e,\alpha}}{1-k}=\epsilon
\end{eqnarray*}
 Hence we have $|| x_m - x_n||^i_{e,\alpha}\leq \epsilon$ and $|\lambda_{e,m} - \lambda_{e,n}| = |\xi_e - \xi_e| =0 < \epsilon$.
 So Definition $\ref{7.8}$ implies that $\{(\tilde x_E)_n\}$ is a F.S Cauchy sequence in $X$ and consequently since $X$ is a F.S Banach space,
  it is a F.S convergent sequence.\\
 Let $(\tilde x_E)_n \rightarrow \tilde x_E$. Therefore Theorems $\ref{th7.8}$ and $\ref{th7.9}$ imply that
$T(\tilde x_E)_n \rightarrow T\tilde x_E$.\\ On the other hand, $(\tilde x_E)_n = T(\tilde x_E)_{n-1} $. So $ T(\tilde x_E)_n \rightarrow \tilde x_E$.\\
 Theorem $\ref{th7.5}$ implies that $\tilde x_E=T\tilde x_E$ i.e., $\tilde x_E$
is a F.S fixed point of $T$ in $X$.\\
Now we show that it is unique. Suppose that  $T\tilde y_E=\tilde y_E$ then Definition $\ref{7.6}$ and Theorem $\ref{th7.9}$ imply that\\
$\bar 0_E \tilde\preceq ||\tilde x_E \tilde\ominus\tilde y_E|| = ||T\tilde x_E \tilde\ominus T\tilde y_E|| \tilde\preceq k ||\tilde x_E \tilde\ominus\tilde y_E||$\\
$\Rightarrow \bar 0_E \tilde\preceq ||\tilde x_E \tilde\ominus\tilde y_E|| \tilde\preceq k ||\tilde x_E \tilde\ominus\tilde y_E||$\\
$\Rightarrow \bar 0_E \tilde\preceq ||\tilde x_E \tilde\ominus\tilde y_E||(1-k) \tilde\preceq \bar 0_E$\\
$\Rightarrow \bar 0_E = ||\tilde x_E \tilde\ominus\tilde y_E||$\\
$\Rightarrow \tilde x_E \tilde\ominus\tilde y_E=\bar 0_E$\\
 It means that for all $e\in E$ and $z\in X$ we have\\
 $(\tilde x \tilde\ominus\tilde y)_e(z)=\left\{\begin{array}{ll}
  \min\{\tilde x_e(x) , \tilde y_e(y)\}=1 & \mbox{if $z=x-y=0$}\\
  0 & \mbox{otherwise}
 \end{array}
 \right.$\\
 So $(\tilde x \tilde\ominus\tilde y)(e)= (\overline{x-y})(e)$\\
$\Rightarrow \tilde x_E \tilde\ominus\tilde y_E={\overline{x-y}}_E=\bar 0_E$\\
  $\Rightarrow x=y$\\
Hence $\tilde x_E$ is the only F.S fixed point of $T$ in $X$. $\Box$
%%%%%%%%%%%%%%%%%%%%%%%%%%%%%%%%%%%%%%%%%%%%%%%%%%%%%%%%%%%%%

\section{ Discussion and Conclusion}\vspace{-4pt}
The language used in our daily life situations is usually full of imprecise phrases. Furthermore, measurement applied in the areas related to economics,
social science, environmental science, and etc., is not crisp and depends on some parameters
  such as our measurement tool, time and place of measurement, individual observer, and etc.
   On the other hand, the collected data from such surveys are not accurate and have
  degree of uncertainty. Theory of fuzzy sets can model these imprecise data while soft set theory provides a useful method to deal with information
    related to some parameters. So fuzzy soft set theory can be used as a framework to approach and model these kind of occasions.\\
Due to determine how much behavior of some objects in an information system are close to each other, the concepts of norm and limit for fuzzy soft sets are
needed.
 To consider the similarity between some geographical regions based on amount of annual rain or behavior of several
 shareholders with respect to some decision parameters related to stock market are examples of application of fuzzy soft norm and
 fuzzy soft convergency in real-life situations.
 In this work, we introduce a norm over the fuzzy soft classes, called F.S norm, to indicate how much the elements of a universal set are
close to each other with regards to some parameters. In fact the fuzzy soft norm may help us to construct the equivalence classes over a set based on
degree of having some parameters.\\
 We firstly introduce the concept of fuzzy soft number for the first time and then consider some arithmetic
operations over them. Then we present fuzzy soft norm and give the relationship between the concept of norm in the sense of classic and F.S norm.
 The concept of F.S sequence and F.S convergent sequence are given. We also introduce a parametrization extension of fixed-point theorem in F.S normed
 spaces.\\
 This paper may be the beginning for future research on fuzzy soft inner product, Banach spaces, fixed-point theorem and etc.
%%%%%%%%%%%%%%%%%%%%%%%%%%%%%%%%%%%%%%%%%%%%%%%%%%%%%%%%%%%%%%%%%%%%%%%%
\section*{Acknowledgement}
The authors acknowledge this research was part of the research project and partially supported by the University Putra Malaysia under the ERGS
1-2013/5527179.


\begin{thebibliography}{11}
\vspace{-7pt}
\bibitem{M}
D. Molodtsov, Soft set theory-first results, Computers and Mathematics with applications 37 (1999) 19-31.

\bibitem{M-B-R}
 P. K. Maji, R. Biswas and A. R. Roy, Soft Set Theory, Computers and Mathematics with applications 45 (2003) 555-562.

\bibitem{M-R}
P. K. Maji and A. R. Roy, An application of soft sets in adecisionmaking problem,Computers and Mathematics with applications 44 (2002) 1077-1083.

\bibitem{m-b-r}
      P. K. Maji,R. Biswas, A. R. Roy, Fuzzy Soft Set ,Journal of Fuzzy Mathematics 9 (3) (2001) 589–602.

\bibitem{m-r}
     A. R. Roy, P. K. Maji, A fuzzy Soft Set theoretic approach to decision making problems, Journal of Computational and Applied Mathematics 203 (2007) 412 – 418

\bibitem{A-C}
     H. Akta\c{s}, N. \c{C}a\v{g}man, Soft sets and soft groups, Information Sciences 177 (2007) 2726-3332.

\bibitem{A-F-M}
   M. I. Ali, F. Feng, X. liu, W.K. Min, On some new operations in soft set theory, Computers and Mathematics with applications 57(9)
(2009) 1547-1553.

\bibitem{S}
 A. Sezgin, A. O. Atagun, On operations of soft sets, Computers and Mathematics with applications
61 (2011) 1457-1467.

\bibitem{R-S}
S. Roy, T. K. Samanta, A note on fuzzy soft topological spaces, Annals of Fuzzy Mathematics and Informatics, Volume 3, No. 2, April 2012, pp. 305- 311,
ISSN 2093–9310.

\bibitem{T-K}
     B. Tanay, M. B. Kandemir, Topological structure of fuzzy soft sets, Computers and Mathematics with Applications 61 (2011) 2952-2957.

\bibitem{M-D}
J. Mahanta  and P. K. Das, Results on fuzzy soft topological spaces, arXiv:1203.0634v1 [math.GM] 3 Mar 2012.

\bibitem{S-Y}
T. Simsekler, S. Yuksel, Fuzzy soft topological spaces, Annals of Fuzzy Mathematics and Informatics, May 2012, ISSN 2093-9310.

\bibitem{A-K}
A. Zahedi Khameneh and A. Kilicman, Separation axioms in Fuzzy Soft topological spaces, submitted for publication.

\bibitem{A-K-S}
     A. Zahedi Khameneh , A. Kilicman and A. R. Salleh, Fuzzy Soft Product Topology, submitted for publication.

\bibitem{Z}
 L.A. Zadeh, Fuzzy Sets, Information and Control 8 (1965) 338-353.

\bibitem{ch}
     C.L. Chang, Fuzzy Topological Spaces, Journal of Mathematical Analysis and Applications 24, 182-190 (1968).

\bibitem{L}
     R. Lowen, Fuzzy Topological Spaces and Fuzzy Compactness, Journal of Mathematical Analysis and Applications 56, 621-633 (1976).

\bibitem{K-S}
     O. Kaleva and S. Seikkala, On Fuzzy Metric Spaces, Fuzzy Sets and Systems 12, 215-229 (1984).

\bibitem{D-H-L}
     J. G. Digkman, H. van Haeringen, and S.J. De Lange, Fuzzy Numbers, Journal of Mathematical Analysis and Applications 92, 301-341 (1983).

\bibitem{M-M}
     M. Itoh, M. Cho, Fuzzy bounded operators, Fuzzy Sets and Systems 93, 353-362 (1998).

\bibitem{L-L-P}
     B-S. Lee, S-J. Lee, K-M. Park, The completions of fuzzy metric spaces and fuzzy normed linear spaces, Fuzzy Sets and Systems 106, 469-473 (1999).

\bibitem{X-Z}
     J. Xiao, X. Zhu, On linearly topological structure and property of fuzzy normed linear space, Fuzzy Sets and Systems 125, 153-161 (2002).

\bibitem{x-z}
     J. Xiao, X. Zhu, Fuzzy normed space of operators and its completeness, Fuzzy Sets and Systems 133, 389-399 (2003).

\bibitem{B-S}
     T. Bag, S.K. Samanta, Fuzzy bounded linear operators in Felbin's type fuzzy normed linear spaces, Fuzzy Sets and Systems 159, 685-707 (2008).

     \bibitem{BA-AKH}
       B. Ahmad and Athar Kharal, On Fuzzy Soft sets, Advances in Fuzzy
     Systems, Volume 2009, Article ID 586507.

\bibitem{kh-a}
     Athar Kharal and  B. Ahmad, Mappings on Fuzzy Soft Classes, Advances in Fuzzy
     Systems, Volume 2009, Article ID 407890.
     3308-3314.

     \bibitem{paw}
     Pawlak. Z, Rough sets, Int J Comput Inform Sci, 177 (1982) 3-27.
     3308-3314.

      \bibitem{pu-liu}
     P. M. PU and Y. M. LIU, Fuzzy topology.I.Neighborhood structure of a fuzzy point and Moore-Smith convergence, J Math. Anal. Appl, 76 (1980) 571-599.
     3308-3314.
\end{thebibliography}
\end{document}